\documentclass[10pt]{article}

\usepackage{amssymb,amsmath,amsthm}
\usepackage{amsfonts}
\usepackage{graphicx}
\usepackage[usenames]{color}
\usepackage{url}
\usepackage{algorithm,algorithmic}
\usepackage{dsfont,verbatim}
\usepackage{enumitem}
\usepackage{epsfig}
\usepackage[margin=1.2in]{geometry}
\graphicspath{{./pics/}}
\allowdisplaybreaks

\newtheorem{theorem}{Theorem}[section]
\newtheorem{lemma}{Lemma}[section]

\numberwithin{equation}{section}

\def\d{{\rm d}}

\def\ii{\textrm{i}}

\title{\bf Convergence of finite element solutions of stochastic partial integro-differential equations driven by white noise} 

\author{Max Gunzburger\thanks{Department of Scientific Computing, Florida State
University, Tallahassee, FL 32306, USA.
(\texttt{gunzburg@fsu.edu})}
\and Buyang Li\thanks{Department of Applied Mathematics, The Hong Kong Polytechnic University, Hung Hom, Kowloon, Hong Kong.
(\texttt{buyang.li@polyu.edu.hk})}
\and Jilu Wang\thanks{Department of Scientific Computing, Florida State University, Tallahassee, FL 32306, USA.(\texttt{jwang13@fsu.edu})}}
\date{}

\begin{document}

\maketitle

\begin{abstract}
Numerical approximation of a stochastic partial integro-differential equation driven by a space-time white noise is studied by truncating a series representation of the noise, with finite element method for spatial discretization and convolution quadrature for time discretization. Sharp-order convergence of the numerical solutions is proved up to a logarithmic factor. 
Numerical examples are provided to support the theoretical analysis. 
\\[5pt]
\end{abstract}


\pagestyle{myheadings}
\thispagestyle{plain}

\section{Introduction}\label{sec:intro}
For given $\alpha\in(0,2)$, we consider the stochastic partial integro-differential equation (PDE)
\begin{align}\label{Frac-SPDE} 
\left\{\begin{aligned}
&\partial_t \psi(x,t)-\Delta \partial_t^{1-\alpha}\psi(x,t) =   
f(x,t) + \sigma \,\dot W(x,t)    
&& (x,t)\in \mathcal{O}\times {\mathbb R}_+ \\
&\partial_t^{1-\alpha}\psi(x,t)=0 && (x,t)\in \partial\mathcal{O}\times {\mathbb R}_+  \\
&\psi(x,0)=\psi_0(x) && x\in \mathcal{O}
\end{aligned}\right.
\end{align}
in a convex polygon/polyhedron $\mathcal{O}\subset\mathbb{R}^d$, $d\in\{1,2,3\}$, 
where $\Delta: H^2(\Omega)\cap H^1_0(\mathcal{O})\rightarrow L^2(\mathcal{O})$ denotes the Laplacian operator, $f(x,t)$ a given deterministic source function, $\psi_0(x)$ a given deterministic initial data, $\sigma$ a given positive parameter, and $\dot W(x,t)$ a space-time white noise, i.e., the time derivative of a cylindrical Wiener process on $L^2(\mathcal{O})$ with an underlying probability sample space $\Omega$. The Caputo fractional time derivative/integral $\partial_t^{1-\alpha}\psi$ is defined by (cf. \cite[pp. 91]{KST}) 
\begin{align} \label{Caputo}
\partial_t^{1-\alpha} \psi(x,t) 
:= \left\{
\begin{aligned}
&\frac{1}{\Gamma(\alpha)}  \int_0^t (t-s)^{\alpha-1} \frac{\partial\psi(x,s)}{\partial s} \d s 
&&\mbox{if}\,\,\, \alpha\in(0,1], \\
&\frac{1}{\Gamma(\alpha-1)} 
\int_0^t (t-s)^{\alpha-2}\psi(x,s)\d s 
&&\mbox{if}\,\,\, \alpha\in(1,2), 
\end{aligned}
\right.
\end{align}
where $\Gamma(s):=\int_0^\infty t^{s-1}e^{-t}\d t$ denotes the Euler Gamma function. 
If $\alpha=1$, then \eqref{Frac-SPDE} recovers the standard stochastic parabolic equation.

Problem \eqref{Frac-SPDE} can be used to describe the behavior of complex  phenomena in mathematical physics, such as viscoelasticity and heat conduction in materials with memory subject to stochastic noises \cite{ClementDaPrato,KP,MijenaNane2015}. For any given initial data $\psi_0\in L^2(\mathcal{O})$ and source $f\in L^1(0,T;L^2(\mathcal{O}))$,  problem \eqref{Frac-SPDE} has a mild solution $\psi\in C([0,T];L^2(\Omega; L^2(\mathcal{O})))$; see Appendix \ref{AppendixA}.

Many efforts have been made in developing efficient numerical methods with rigorous error analyses for solving \eqref{Frac-SPDE}, with or without the stochastic noise. In \cite{LubichSloanThomee:1996}, Lubich et.al. have considered the deterministic version of this problem in the case $\alpha\in(1,2)$. The discretization used convolution quadrature (CQ) based on backward difference methods (BDFs) in time and piecewise linear finite elements in space. The authors have proved optimal-order convergence rate of the numerical scheme for nonsmooth initial data. To achieve higher-order temporal convergence rates, the CQ generated by second-order BDF and Crank-Nicolson methods have been considered in \cite{CuestaLubichPalencia:2006} and \cite{JinLiZhou-CN} for solving \eqref{Frac-SPDE} and its equivalent formulation, respectively. Due to the singularity of the solution of fractional evolution PDEs, the standard BDF and Crank-Nicolson CQs need to be corrected at several initial steps to achieve the desired order of convergence. Initial correction of higher-order BDF methods for fractional evolution PDEs has been considered in \cite{JinLiZhou-BDF} recently. Compared with the deterministic problem, the major technical difficulties in the development and analyses of numerical schemes for \eqref{Frac-SPDE} are due to the space-time white-noise forcing, which leads to low regularity of the solution in both time and space. 
In the case $\alpha=1$, Allen et. al. \cite{AllenNovoselZhang1998} developed a fully discrete numerical scheme for solving stochastic parabolic problem, for which the white noise was approximated by piecewise constant random processes and a sharp order of convergence was proved. See also Du and Zhang \cite{DuZhang2002} for some special noises, Shardlow \cite{Shardlow1999} for the space-time white noise discretized by the spectral method, and Yan \cite{Yan} for a nonlinear stochastic parabolic system with Wiener process discretized by the generalized $L^2$-projection operator. In \cite{KP}, Kov\'acs and Printems developed a CQ based on backward Euler method for the model \eqref{Frac-SPDE} with $\alpha\in(1,2)$,  where the $Q$-Wiener process was discretized by the generalized $L^2$-projection operator. For the space-time white noise case, a strong order of convergence of the numerical solution was proved in one-dimensional spatial domains, i.e.,  
\begin{align*}
\mathbb{E}\|\psi(\cdot,t_n)-\psi_n^{(h)}\|_{L^{2}(\mathcal{O})}
=O(\tau^{\frac{1}{2}-\frac{\alpha}{4}-\varepsilon}+h^{\frac{1}{\alpha}-\frac{1}{2}-\varepsilon})
\quad\mbox{for}\,\,\, \alpha\in(1,2)\,\,\,\mbox{and}\,\,\, d=1,
\end{align*}
where $\varepsilon$ can be arbitrarily small, $\psi(\cdot,t_n)$ and $\psi_n^{(h)}$ denote the PDE's mild solution and numerical solution at time $t_n$, respectively, 
$\tau$ denotes the temporal step size, and $h$ denotes the spatial mesh size.  
For $\alpha\in(0,2)$, a sharp order convergence rate $O(\tau^{\frac{1}{2}-\frac{\alpha d}{4}})$ was proved in \cite{GunzburgerLiWang2017} for a CQ time discretization of \eqref{Frac-SPDE} in general $d$-dimensional spatial domains, with $d\in\{1,2,3\}$, without the deterministic forcing. We refer the readers to \cite{AntonCohenLarssonWang2016,BanasBrzeznikProhl2013,FengLiProhl2014} for numerical analysis of other nonlinear physical stochastic equations.

This article is a continuation of \cite{GunzburgerLiWang2017} in the spatially discrete setting, by truncating a series representation of the space-time white noise and solving the truncated problem by the finite element method. For the resulting fully discrete numerical scheme, we prove the sharp-order convergence 
\begin{align}
\mathbb{E}\|\psi(\cdot,t_n)-\psi_n^{(h)}\|_{L^{2}(\mathcal{O})}
=
\left\{\begin{aligned}
&O\big(\tau^{\frac{1}{2}-\frac{\alpha d}{4}}
+\ell_h^{\frac{1}{2}} h^{\frac{1}{\alpha}-\frac{d}{2}}\big) 
&&\mbox{if}\,\,\, \alpha\in\Big[\frac{1}{2},\frac{2}{d}\Big) , \\[5pt]
&O\big(\tau^{\frac{1}{2}-\frac{\alpha d}{4}}+ h^{2-\frac{d}{2}} \big) 
&&\mbox{if}\,\,\, \alpha\in\Big(0,\frac{1}{2}\Big) ,
\end{aligned}\right.  
\end{align}
up to a logarithmic factor $\ell_h^{\frac{1}{2}}=(\ln(e+1/h))^{\frac{1}{2}}$, in general $d$-dimensional spatial domains, $d\in\{1,2,3\}$. 
The main contributions of this paper are the following. 
\begin{enumerate}[label={\rm(\arabic*)},ref=\arabic*,topsep=2pt,itemsep=0pt,partopsep=1pt,parsep=1ex,leftmargin=20pt]
\item
Sharper-order spatial convergence is proved in the case $\alpha\in(1,2)$ and $d=1$ (up to a logarithmic factor $\ell_h^{\frac{1}{2}}$). 

\item The error estimates are extended to $\alpha\in(0,\frac{2}{d})$ and multi-dimensional domains. 

\item An interesting phenomenon is found: 
the spatial order of convergence $\frac{1}{\alpha}-\frac{d}{2}$ increases to ${2-\frac{d}{2}}$ as $\alpha$ decreases to $\frac{1}{2}$, 
and stays at this order when $\alpha$ further decreases.  

\item Less regularity assumption on $f$: 
the error estimates in the literature all rely on certain regularity of $\frac{\partial f}{\partial t}$ (even for the deterministic problems, cf. \cite[Theorem 3.6]{JinLazarovZhouSISC2016} and \cite[Theorem 3.3]{LubichSloanThomee:1996}). We relax such conditions to an optimal integrability condition $f\in L^{\frac{4}{2+\alpha d},1}(0,T;L^2(\mathcal{O}))$ to match the convergence rate of the stochastic problem. 
Consequently, the source $f$ does not need to be continuous in time. 
\end{enumerate}

The rest of this paper is organized as follows. In section 2, we recall some basic preliminary results, introduce the numerical scheme for problem \eqref{Frac-SPDE}, and state the main results. Based on an integral representation of the numerical solution and careful analyses of the resolvent operator, the strong convergence rates are proved in section 3 and section 4. 
Numerical examples are given in section 5 to illustrate the theoretical results. 

Throughout this paper, we denote by $C$, with or without a subscript/superscript, a generic constant independent of $n$, $\tau$, and $h$, which could be different at different occurrences.
%

\section{The main results}\label{sec:main}

\subsection{Notations}\label{sec:notation}
We denote by $(\cdot,\cdot)$ and $\|\cdot\|$ the inner product and norm of $L^2(\mathcal{O})$, respectively. The operator norm on $L^2(\mathcal{O})$ is also denoted by $\|\cdot\|$ (as it is induced by the norm of $L^2(\mathcal{O})$).
Let $\dot H^s(\mathcal{O})\subset L^2(\mathcal{O})$ denote the Hilbert space induced by the norm 
\begin{align}\label{Def-frac-norm}
\|\varphi\|_{\dot H^s(\mathcal{O})}
:=\sum_{j=1}^\infty \lambda_j^{2s}|(\varphi,\phi_j)|^2 ,
\end{align}
where $\phi_j$, $j=1,2,\dots$, denote the $L^2$-norm normalized eigenfunctions of the Laplacian operator $-\Delta$ corresponding to the eigenvalues $\lambda_j$ , $j= 1, 2,\dots$, arranged in nondecreasing order. In particular, $\dot H^0(\mathcal{O})=L^2(\mathcal{O})$, 
$\dot H^1(\mathcal{O})=H^1_0(\mathcal{O})$ and $\dot H^2(\mathcal{O})=H^2(\mathcal{O})\cap H^1_0(\mathcal{O})$; see \cite{Thomee:2006}. 
For $1<p<\infty$ we denote by $L^{p,1}(0,T;L^2(\mathcal{O}))$ the standard Lorentz space of functions 
defined on $\mathcal{O}\times(0,T)$ (see \cite[section 1.4]{Grafakos2008}), satisfying  
\begin{align}\label{Lp1-use}
\sup_{t\in(0,T)}\int_0^t \!\!(t-s)^{-\frac{1}{p'}} \|f(\cdot,s)\| \d s 
\le C\|f\|_{L^{p,1}(0,T;L^2(\mathcal{O}))}  \quad\forall\, f\in L^{p,1}(0,T;L^2(\mathcal{O})), 
\end{align}
where $p'$ denotes the dual of $p$, i.e., $\frac{1}{p'}+\frac{1}{p}=1$. 
For $1<p<\infty$, the Lorentz space $L^{p,1}(0,T;L^2(\mathcal{O}))$ is an intermediate real interpolation space between $L^{1}(0,T;L^2(\mathcal{O}))$ and $L^{\infty}(0,T;L^2(\mathcal{O}))$ (see \cite[Theorem 5.2.1]{BL1976}), satisfying 
\begin{align}
L^{q}(0,T;L^2(\mathcal{O}))\hookrightarrow L^{p,1}(0,T;L^2(\mathcal{O})) \quad \forall\, q>p\ge 1 .
\end{align}

Let $\{t_n=n\tau\}_{n=0}^N$ denote a uniform partition of the time interval $[0,T]$, 
with a step size $\tau=T/N$, and $u^n=u(x,t_n)$. 
If we denote by $f_\tau$ the following function (piecewise constant in time):
\begin{align}\label{f-tau}
f_\tau (\cdot,s)= \frac{1}{\tau}\int_{t_{n-1}}^{t_n}f(\cdot, t)\d t \quad \forall\, s\in (t_{n-1},t_n] , \,\,\, n=1,2,\dots,N ,
\end{align}  
then it is well known that 
\begin{align}
&\|f_\tau\|_{L^{1}(0,T;L^2(\mathcal{O}))}\le \|f\|_{L^{1}(0,T;L^2(\mathcal{O}))}
&& \forall\, f\in L^{1}(0,T;L^2(\mathcal{O})), \\
&\|f_\tau\|_{L^{\infty}(0,T;L^2(\mathcal{O}))}\le \|f\|_{L^{\infty}(0,T;L^2(\mathcal{O}))}
&& \forall\, f\in L^{\infty}(0,T;L^2(\mathcal{O})) .
\end{align}  
The real interpolation of the last two inequalities yields (see \cite[Definition 2.4.1]{BL1976} and \cite[Theorem 5.2.1]{BL1976})
\begin{align}\label{Lp1-stable}
&\|f_\tau\|_{L^{p,1}(0,T;L^2(\mathcal{O}))}\le C\|f\|_{L^{p,1}(0,T;L^2(\mathcal{O}))} 
&& \forall\, f\in L^{p,1}(0,T;L^2(\mathcal{O})) .
\end{align}  
The last inequality will be used in this paper.

For $\alpha\in(0,1]$, we approximate the Caputo fractional time derivative $\partial_t^{1-\alpha} u(x,t_n)=\partial_t^{1-\alpha} (u(x,t_n)-u(x,0))$ by the backward Euler CQ (cf. \cite[(2.4)]{JinLiZhou-CN} and \cite{Lubich:1988-1,Lubich:1988-2,LubichSloanThomee:1996,Sanz-Serna1988}):  
\begin{align}\label{convqdr}
\bar\partial_\tau^{1-\alpha} (u_n-u_0)
:=\frac{1}{\tau^{1-\alpha} } \sum_{j=1}^n b_{n-j}(u_j-u_0) ,\quad n=1,2,\dots,N.
\end{align}
For $\alpha\in(1,2)$, we approximate the Caputo fractional time derivative $\partial_t^{1-\alpha} u(x,t_n)$ by the CQ without subtracting the initial data (cf. \cite[(1.15)]{LubichSloanThomee:1996}), i.e., 
\begin{align}\label{convqdr-2}
\bar\partial_\tau^{1-\alpha}u_n 
:=\frac{1}{\tau^{1-\alpha} } \sum_{j=1}^n b_{n-j} u_j  ,\quad n=1,2,\dots,N .
\end{align}
In both \eqref{convqdr} and \eqref{convqdr-2}, the coefficients $b_j$, $j=0,1,2,\dots$, are determined by the power series expansion 
$$
(1-\zeta)^{1-\alpha}=\sum_{j=0}^\infty b_{j}\zeta^j  \quad\forall\,|\zeta|<1 ,\,\, \zeta\in{\mathbb C}.
$$ 
Besides, we define the standard backward Euler difference operator
\begin{align}\label{def-partial_tau}
\bar\partial_\tau u_n
:=\frac{u_n-u_{n-1}}{\tau} ,\quad n=1,2,\dots,N .
\end{align} 
The complex-valued function 
\begin{align}
\label{definition-d}
\delta_\tau(\zeta)= \frac{1-\zeta}{\tau}  \quad \textrm{for }\zeta\in {\mathbb C}\backslash[1,\infty)  
\end{align}
is called the generating function of the backward Euler difference operator. It plays an important role in the analysis of the CQ. In particular, for any sequence $\{v_n\}_{n=0}^\infty\in \ell^2(L^2(\mathcal{O}))$ we have 
\begin{equation}\label{generate-dv}
\begin{aligned}
\sum_{n=1}^\infty (\bar\partial_\tau^{1-\alpha} v_n)\zeta^n
&=\sum_{n=1}^\infty  \frac{1}{\tau^{1-\alpha}} \sum_{j=1}^n b_{n-j}v_j\zeta^n
=(\delta_\tau(\zeta))^{1-\alpha}\sum_{j=1}^\infty v_j\zeta^j ,\quad\forall\,|\zeta|<1 .
\end{aligned} 
\end{equation}

Let $\mathcal{T}_h$ be a quasi-uniform triangulation of the domain $\mathcal{O}$ into $d$-dimensional simplexes $\pi_h$, $\pi_h\in\mathcal{T}_h$, with a mesh size $h$ such that $0<h<h_0$ for some constant $h_0$. A continuous piecewise linear finite element space $X_h$ over the triangulation $\mathcal{T}_h$ is defined by
$$
X_h=\{\phi_h\in H_0^1(\mathcal{O}): \phi_h|_{\pi_h}\mbox{ is a linear function}, \forall \pi_h\in\mathcal{T}_h \}.
$$
Over the finite element space $X_h$, we denote the $L^2$ projection $P_h: L^2(\mathcal{O})\rightarrow X_h$ and Ritz projection $R_h:H_0^1(\mathcal{O})\rightarrow X_h$ by 
\begin{align*}
\begin{aligned}
&(P_h\varphi,\phi_h)=(\varphi,\phi_h) && \forall \phi_h\in X_h,  \\
&(\nabla R_h\varphi,\nabla\phi_h)=(\nabla\varphi,\nabla\phi_h) && \forall \phi_h\in X_h .
\end{aligned}
\end{align*}

It is well known that the $L^2$ projection and Ritz projection satisfy the following standard error estimates (\cite{Thomee:2006}):
\begin{align}\label{StabPh}
&\|P_h\phi\|  \le C\|\phi\|  &&\forall\,\phi\in L^2(\mathcal{O}), \\
&\|P_h\phi-\phi\|  \le Ch^\gamma\|\phi\|_{\dot H^\gamma(\mathcal{O})}  &&\forall\,\phi\in \dot H^\gamma(\mathcal{O}),\,\,\,\gamma\in[0,2], \label{StabPh-} \\
&\|P_h\phi-R_h\phi\|  \le Ch^2\|\phi\|_{\dot H^2(\mathcal{O})}   &&\forall\,\phi\in \dot H^2(\mathcal{O}) . 
\label{StabPh2}
\end{align}
Through defining the discrete Laplacian $\Delta_h:X_h\rightarrow X_h$ by 
$$
(\Delta_h\varphi_h,\phi_h)=-(\nabla\varphi_h,\nabla\phi_h) 
\quad \forall\varphi_h,\phi_h\in X_h
$$
and using the inverse inequality, the inequality \eqref{StabPh2} implies 
\begin{align}
&\|\Delta_h(P_h\phi-R_h\phi)\| 
\le Ch^{-2} \|P_h\phi-R_h\phi\| \le C\|\phi\|_{\dot H^2(\mathcal{O})}  &&\forall\,\phi\in \dot H^2(\mathcal{O}) .
\label{StabPh3}
\end{align}
Since $\Delta_hR_h\phi=P_h\Delta\phi$, it follows that  
\begin{align}\label{DeltahPh}
\|\Delta_h P_h\phi\| 
&\le  \|\Delta_hR_h\phi \| + \|\Delta_h(P_h\phi-R_h\phi)\| \nonumber\\
&= \|P_h\Delta\phi \| + \|\Delta_h(P_h\phi-R_h\phi)\| \nonumber\\
&\le C\|\phi\|_{\dot H^2(\mathcal{O})}  &&\forall\,\phi\in \dot H^2(\mathcal{O}) . 
\end{align}
The complex interpolation between \eqref{StabPh} and \eqref{DeltahPh} yields 
\begin{align}\label{StabPh5}
\|\Delta_h^\gamma P_h\phi\| 
&\le C\|\phi\|_{\dot H^{2\gamma}(\mathcal{O})}  &&\forall\,\phi\in \dot H^{2\gamma}(\mathcal{O}) ,\,\,\,\gamma\in[0,1].  
\end{align}
Similarly, the complex interpolation between \eqref{StabPh2} and \eqref{StabPh3} yields 
\begin{align}\label{StabPh4}
\|\Delta_h^\gamma (P_h\phi-R_h\phi)\| 
&\le Ch^{2-2\gamma}\|\phi\|_{\dot H^2(\mathcal{O})}  &&\forall\,\phi\in \dot H^{2}(\mathcal{O}) ,\,\,\,\gamma\in[0,1].  
\end{align}
The estimates \eqref{StabPh}-\eqref{StabPh4} will be frequently used in this paper.


\subsection{The numerical scheme and main theorem}

Recall that the cylindrical Wiener process on $L^2(\mathcal{O})$ can be represented as 
(cf. \cite[Proposition 4.7, with $Q=I$ and $U_1$ denoting some negative-order Sobolev space]{PratoZabczyk2014})
\begin{align*}
W(x,t) =\sum_{j=1}^\infty \phi_j(x) W_j(t)  
\end{align*} 
with independent one-dimensional Wiener processes $W_j(t)$, $j=1,2,\dots$. 
We approximate the space-time white noise $\dot W(x,t)$ by 
\begin{align*}
\bar\partial_\tau W^M(x,t_n)
=\sum_{j=1}^M \phi_j(x) \bar\partial_\tau W_j(t_n)
\end{align*}
with $M:=[h^{-d}]+1$, the largest integer that does not exceed $h^{-d}+1$.  
Clearly, we have 
\begin{align}
h^{-d}\le M\le Ch^{-d} \quad \forall\, 0<h<h_0 ,
\end{align}
where the constant $C$ may depend on $h_0$.

With the above notations, we propose the following fully discrete scheme for problem \eqref{Frac-SPDE}: find $$
\psi_n^{(h)}\in
\left\{\begin{aligned}
&\psi_0^{(h)}+ X_h &&\mbox{in the case $\alpha\in(0,1]$}\\
&X_h &&\mbox{in the case $\alpha\in(1,2)$}
\end{aligned}\right.
\qquad
n=1,2,\dots,N,
$$
with $\psi_0^{(h)}=P_h\psi_0$, 
such that the following equations are satisfied for all $\phi_h\in X_h$: 
\begin{align}\label{fully-discrete}
&\left(\bar\partial_\tau \psi_n^{(h)}  ,\phi_h\right) +\left(\nabla\bar\partial_\tau^{1-\alpha} (\psi_n^{(h)}-\psi_0^{(h)}) ,\nabla\phi_h\right) \nonumber \\
&\,\,=
\left(f_n,\phi_h\right) 
+\left(\sigma \bar\partial_\tau W^M(\cdot,t_n),\phi_h \right), && 
\mbox{if $\alpha\in(0,1]$}  ,\\[15pt]
&\left(\bar\partial_\tau \psi_n^{(h)}  ,\phi_h\right) +\left(\nabla\bar\partial_\tau^{1-\alpha} \psi_n^{(h)} ,\nabla\phi_h\right) \nonumber \\
&\,\,=
\left(f_n,\phi_h\right) 
+\left(\sigma \bar\partial_\tau W^M(\cdot,t_n),\phi_h \right), &&
\mbox{if $\alpha\in(1,2)$}   ,
\label{fully-discrete-2}
\end{align}
where $f_n$ is the average of $f$ over the subinterval $(t_{n-1},t_n]$, i.e., 
\begin{align}\label{Def-fn} 
f_n=\frac{1}{\tau}\int_{t_{n-1}}^{t_n}f(\cdot,t) \d t . 
\end{align}

Through the discrete Laplacian $\Delta_h$, 
we can rewrite the fully discrete scheme \eqref{fully-discrete}-\eqref{fully-discrete-2} in the following equivalent forms:
\begin{align}\label{fully-discrete2} 
\begin{aligned} 
&\bar\partial_\tau \psi_n^{(h)} -\Delta_h \bar\partial_\tau^{1-\alpha}(\psi_n^{(h)} -\psi_0^{(h)} ) 
= P_hf_n+ \sigma P_h\bar\partial_\tau W^M(\cdot,t_n) , &&\mbox{if}\,\,\,\alpha\in(0,1],\\
&\bar\partial_\tau \psi_n^{(h)}  -\Delta_h \bar\partial_\tau^{1-\alpha} \psi_n^{(h)}
= P_hf_n+\sigma  P_h\bar\partial_\tau W^M(\cdot,t_n) , &&\mbox{if}\,\,\,\alpha\in(1,2).
\end{aligned}
\end{align}
Note that $\psi_0^{(h)}\in L^2(\Omega;X_h)$ and $P_hf_n+\sigma  P_h\bar\partial_\tau W^M(\cdot,t_n)\in L^2(\Omega;X_h)$ for $n=1,2,\dots,N$. If the numerical solutions $\psi_n^{(h)}\in L^2(\Omega;X_h)$, $n=0,1,\dots,m-1$, then we define 
$$g_m^{(h)}
:=
\left\{\begin{aligned}
&\tau^{\alpha-1}\sum_{j=1}^{m-1}b_{m-j}\Delta_h(\psi_j^{(h)}-\psi_0^{(h)})+P_hf_m+\sigma  P_h\bar\partial_\tau W^M(\cdot,t_m) &&\mbox{if}\,\,\,\alpha\in(0,1],\\
&\tau^{\alpha-1}\sum_{j=1}^{m-1}b_{m-j}\Delta_h\psi_j^{(h)}+P_hf_m+\sigma  P_h\bar\partial_\tau W^M(\cdot,t_m) &&\mbox{if}\,\,\,\alpha\in(1,2) .
\end{aligned}
\right.
$$
Then $g_m^{(h)}\in L^2(\Omega;X_h)$, and the numerical solution defined by \eqref{fully-discrete2} is given by 
$$
\psi_m^{(h)}
=
\left\{\begin{aligned} 
&\psi_0^{(h)}+(\tau^{-1}-\tau^{\alpha-1}b_0\Delta_h)^{-1}\big(\tau^{-1}(\psi_{m-1}^{(h)}-\psi_0^{(h)})+g_m^{(h)} \big) &&\mbox{if}\,\,\,\alpha\in(0,1) ,\\
&(\tau^{-1}-\tau^{\alpha-1}b_0\Delta_h)^{-1}\big(\tau^{-1}\psi_{m-1}^{(h)} + g_m^{(h)} \big) &&\mbox{if}\,\,\,\alpha\in(1,2) ,
\end{aligned} 
\right.
$$ 
which is well defined in $L^2(\Omega;X_h)$. By induction, the numerical solutions $\psi_n^{(h)}\in L^2(\Omega;X_h)$, $n=1,2,\dots,N$, are well defined.

The main result of this paper is the following theorem. 
\begin{theorem}\label{MainTHM}
Let $\alpha\in(0,\frac{2}{d})$ with $d\in\{1,2,3\}$, 
$f\in L^{\frac{4}{2+\alpha d},1}(0,T;L^2(\mathcal{O}))$ and 
$\psi_0\in \dot H^{\chi}\!(\mathcal{O})$, with the notations
\begin{align}
&\chi=\min\left(2-\frac{d}{2} \, ,\, \frac{1}{\alpha}-\frac{d}{2}\right)  ,  \label{Def-chi}\\
&\Upsilon(\psi_0,f) = \|\psi_0\|_{\dot H^{\chi}\!(\mathcal{O})}
\!\! + \! \|f\|_{L^{\frac{4}{2+\alpha d},1}(0,T;L^2(\mathcal{O}))} , 
\end{align}
and assume that the spatial mesh size satisfies $0<h<h_0$ for some constant $h_0$. 
Then the numerical solution given by \eqref{fully-discrete2} converges to the mild solution of \eqref{Frac-SPDE} with sharp order of convergence, i.e., 
\begin{align*}
&\mathbb{E}\|\psi(\cdot,t_n)-\psi_n^{(h)}\|_{L^{2} (\mathcal{O})}  
\! \le \! 
\left\{\begin{aligned}
&\!C\big(\sigma\!+\! \Upsilon(\psi_0,f) \big) \big(\tau^{\frac{1}{2}-\frac{\alpha d}{4}}
\! + \! \ell_h^{\frac{1}{2}} h^{\frac{1}{\alpha}-\frac{d}{2}}\big) 
&&\!\!\mbox{if}\ \alpha\in\Big[\frac{1}{2},\frac{2}{d}\Big) ,\\
&\!C\big(\sigma \!+\! \Upsilon(\psi_0,f) \big)\big(\tau^{\frac{1}{2}-\frac{\alpha d}{4}} \! + h^{2-\frac{d}{2}} \big) 
&&\!\!\mbox{if}\ \alpha\in\Big(0,\frac{1}{2}\Big) , 
\end{aligned}\right.  
\end{align*}
where $\mathbb{E}$ denotes the expectation operator, $\ell_h=\ln(e+1/h)$, 
the constant $C$ is independent of $h$, $\tau$, $n$, $\psi_0$, and $f$ (but may depend on $T$ and $h_0$). 
\end{theorem}

{\bf Proof.} 
Without loss of generality, we can assume $\sigma=1$ in the proof of Theorem \ref{MainTHM}. 
The solution of \eqref{Frac-SPDE} can be decomposed into the solution of the deterministic problem 
\begin{align}\label{Deter-SPDE2} 
\left\{
\begin{array}{ll}
\partial_t v(x,t) -\Delta \partial_t^{1-\alpha}v(x,t) = f (x,t)   &  (x,t)\in \mathcal{O}\times {\mathbb R}_+   \\
\partial_t^{1-\alpha}v(x,t)=0  & (x,t)\in \partial\mathcal{O}\times {\mathbb R}_+     \\
v(x,0)=\psi_0(x)   & x\in \mathcal{O}     
\end{array}\right.
\end{align}
plus the solution of the stochastic problem
\begin{align}\label{Frac-SPDE2} 
\left\{\begin{aligned}
&\partial_t u(x,t) -\Delta \partial_t^{1-\alpha} u(x,t) = \dot W(x,t) 
&& (x,t)\in \mathcal{O}\times {\mathbb R}_+ \\
&\partial_t^{1-\alpha}u(x,t)=0    &&(x,t)\in \partial\mathcal{O}\times {\mathbb R}_+ \\ 
&u(x,0)=0   && x\in \mathcal{O}  . \\ 
\end{aligned}\right.
\end{align} 
Similarly, the solution of \eqref{fully-discrete2} can be decomposed into the solution of the deterministic finite element equation
\begin{equation}\label{Deter-discreteh}
\left\{\begin{aligned}
&\bar\partial_\tau v_n^{(h)}  -\Delta_h \bar\partial_\tau^{1-\alpha}  (v_n^{(h)}-v_0^{(h)}) 
= P_hf_n \\
&v_0^{(h)}=P_h\psi_0 
\end{aligned}
\right.
\qquad\,\, \mbox{if\, $\alpha\in(0,1]$},
\end{equation}
or 
\begin{equation}\label{Deter-discreteh-2}
\left\{\begin{aligned}
&\bar\partial_\tau v_n^{(h)}  -\Delta_h \bar\partial_\tau^{1-\alpha} v_n^{(h)}
= P_hf_n \\
&v_0^{(h)}=P_h\psi_0 
\end{aligned}
\right.\qquad\qquad\qquad \mbox{if\, $\alpha\in(1,2)$},
\end{equation}
plus the solution of the stochastic finite element equation
\begin{equation}\label{Frac-discreteh}
\left\{\begin{aligned}
&\bar\partial_\tau u_n^{(h)}  -\Delta_h \bar\partial_\tau^{1-\alpha} u_n^{(h)}
=  P_h\bar\partial_\tau W^M(\cdot,t_n) \\
&u_0^{(h)}= 0 .
\end{aligned}
\right.
\end{equation}
In the next two sections, we prove Theorem \ref{MainTHM} by estimating 
$\mathbb{E}\|u(\cdot,t_n)-u_n^{(h)}\|$ and $\|v(\cdot,t_n)-v_n^{(h)}\|$, separately. 
In particular, Theorem \ref{MainTHM} follows from \eqref{error-u-unh} and \eqref{error-v-vnh} (in \eqref{error-v-vnh}, we have $\chi=\frac{1}{\alpha}-\frac{d}{2}$ for $\alpha\in[\frac12,\frac2d)$ and $\chi= 2-\frac{d}{2}$ for $\alpha\in(0,\frac12)$).

\section{Stochastic problem: estimate of $\mathbb{E}\|u(\cdot,t_n)-u_n^{(h)}\|$} In this section, we prove the following error estimate for the solutions of \eqref{Frac-SPDE2} and \eqref{Frac-discreteh}:
\begin{align}\label{error-u-unh}
\mathbb{E}\|u(\cdot,t_n)-u_n^{(h)}\|
&\le
\left\{\begin{aligned}
&C(\tau^{\frac{1}{2}-\frac{\alpha d}{4}} 
+ \ell_h^{\frac{1}{2}} h^{\frac{1}{\alpha}-\frac{d}{2}})  &&\mbox{if}\,\,\, \alpha\in\Big[\frac{1}{2},\frac{2}{d}\Big) , \\ 
&C(\tau^{\frac{1}{2}-\frac{\alpha d}{4}} 
+  h^{2-\frac{d}{2}})&&\mbox{if}\,\,\, \alpha\in\Big(0,\frac{1}{2}\Big) .
\end{aligned}\right.
\end{align}
To this end, we introduce a time-discrete system of PDEs:   
\begin{equation}\label{CQ-scheme2}
\left\{
\begin{aligned}
&\bar\partial_\tau u_n -\Delta \bar\partial_\tau^{1-\alpha} u_n 
= \bar\partial_\tau W(\cdot,t_n) ,\quad n=1,2,\dots,N,\\
&u_0=0.
\end{aligned}
\right.
\end{equation}
Then \eqref{Frac-discreteh} can be viewed as the spatially finite element discretization of \eqref{CQ-scheme2},  
and the error can be decomposed into two parts: 
\begin{align}\label{error-split}
\mathbb{E}\|u(\cdot,t_n)-u_n^{(h)}\|
&\le
\mathbb{E}\|u(\cdot,t_n)-u_n\|
+\mathbb{E}\|u_n-u_n^{(h)}\|, 
\end{align}
where the first part on the right-hand side has been estimated in \cite{GunzburgerLiWang2017} (in \cite{GunzburgerLiWang2017} we have only considered zero initial condition $u(\cdot,0)=0$, and in this case the boundary condition $u=0$ on $\partial\Omega$ is equivalent to $\partial_t^{1-\alpha}u=0$ on $\partial\Omega$), i.e.,  
\begin{align}\label{error-time}
\mathbb{E}\|u(\cdot,t_n)-u_n\|
\le
C\tau^{\frac{1}{2}-\frac{\alpha d}{4}} 
\qquad\forall\, \alpha\in\big(0,2/d\big),\,\,\, d\in\{1,2,3\} .
\end{align}
It remains to prove the following estimate in the next three subsections: 
\begin{align}\label{space-discr}
\mathbb{E}\|u_n-u_n^{(h)}\|
\le
\left\{\begin{aligned}
&C\ell_h^{\frac{1}{2}} h^{\frac{1}{\alpha}-\frac{d}{2}}  &&\mbox{if}\,\,\, \alpha\in\Big[\frac{1}{2},\frac{2}{d}\Big) , \\ 
&C h^{2-\frac{d}{2}} &&\mbox{if}\,\,\, \alpha\in\Big(0,\frac{1}{2}\Big) .
\end{aligned}\right.
\end{align}

\subsection{Integral representations}
We estimate $\mathbb{E}\|u_n-u_n^{(h)}\|$ by using integral representations of $u_n$ and $u_n^{(h)}$, respectively. 
We first introduce some notations:  
\begin{align}
& 
\Gamma_{\theta,\kappa} 
\, =\left\{z\in \mathbb{C}: |z|=\kappa  , |\arg z|\le \theta\right\}\cup
  \{z\in \mathbb{C}: z=\rho e^{\pm {\rm i}\theta}, \rho\ge \kappa 
  \}  , \\
&
\Gamma_{\theta,\kappa}^{(\tau)}
\, =\left\{z\in \Gamma_{\theta,\kappa} : |{\rm Im}(z)|\le \frac{\pi}{\tau}\right\} ,
\label{Gamma-tau}
\end{align}
which are contours on the complex plane, oriented with increasing imaginary parts. 
On the truncated contour $\Gamma_{\theta,\kappa}^{(\tau)}$, 
 the following estimates hold. 
\begin{lemma}[\cite{GunzburgerLiWang2017}] \label{ineq-1} 
Let $\alpha\in\big(0,\frac{2}{d}\big)$, $\theta\in\big(\frac{\pi}{2},{\rm arccot}(-\frac{2}{\pi})\big)$ and $\kappa=\frac{1}{T}$ be given, with $\delta_\tau(\zeta)$ defined in \eqref{definition-d}. 
Then 
\begin{align}
&\delta_\tau(e^{-z\tau})\in\Sigma_{\theta} &&\forall\, z\in\Gamma_{\theta,\kappa}^{(\tau)} 
\label{angle-delta} \\
&C_0|z|\le |\delta_\tau(e^{-z\tau})|\le C_1|z|
&&\forall\, z\in\Gamma_{\theta,\kappa}^{(\tau)}  \label{z-delta}  \\
&|\delta_\tau(e^{-z\tau})-z|\le C\tau|z|^2
&&\forall\, z\in\Gamma_{\theta,\kappa}^{(\tau)} \\
&|\delta_\tau(e^{-z\tau})^\alpha-z^\alpha|\le C\tau|z|^{\alpha+1} 
&&\forall\, z\in\Gamma_{\theta,\kappa}^{(\tau)},  \label{zalpha-delta} 
\end{align}
where $\Sigma_\theta:=\{z\in\mathbb{C}\backslash \{0\}:|\arg z|\le\theta<\pi \}$, the constants $C_0,C_1$ and $C$ are independent of $\tau$ and $\kappa$. 
\end{lemma}

Let $\bar\partial_\tau W$ denote a piecewise constant function in time, defined by 
\begin{align*}
&\bar\partial_\tau W(\cdot,t_0):=0 \\
&\bar\partial_\tau W(\cdot,t):=\frac{W(\cdot,t_n)-W(\cdot,t_{n-1})}{\tau}&&\mbox{for}\,\,\, t\in(t_{n-1},t_n],\,\,\, n=1,2,\dots,N.
\end{align*}

Similarly, we define 
\begin{align*}
&\bar\partial_\tau W_j(t_0):=0 \\
&\bar\partial_\tau W_j(t):=\frac{W_j(t_n)-W_j(t_{n-1})}{\tau}&&\mbox{for}\,\,\, t\in(t_{n-1},t_n],\,\,\, n=1,2,\dots,N. 
\end{align*}

Then the following results hold. 
\begin{lemma}\label{un-rep}
Let $\alpha\in\big(0,\frac{2}{d}\big)$ and $\delta_\tau(\zeta)$ be defined as in \eqref{definition-d} with the parameters $\kappa$ and $\theta$ satisfying the conditions of Lemma \ref{ineq-1}.  
\begin{enumerate}
[label={\rm(\arabic*)},ref=\arabic*]\itemsep=5pt
\item 
The solution of the time-discrete problem \eqref{CQ-scheme2} can be represented by 
\begin{align}\label{u-semi}
u_n  
&=
\int_0^{t_n}  
E_\tau(t_n-s)\bar\partial_\tau W(\cdot,s)\d s 
= \sum_{j=1}^\infty \int_{0}^{t_n} E_\tau(t_n-s) \phi_j \bar\partial_\tau W_j(s) \d s  ,
\end{align} 
where the operator $E_\tau(t)$ is given by 
\begin{equation}\label{eqn:EF2}
E_\tau(t) \phi:=\frac{1}{2\pi {\rm i}}\int_{\Gamma_{\theta,\kappa}^{(\tau)}}e^{zt} \frac{z\tau }{e^{z\tau}-1}\delta_\tau(e^{-z\tau})^{\alpha-1}(\delta_\tau(e^{-z\tau})^\alpha  -\Delta )^{-1}\phi\, \d z 
\end{equation}
for $\phi\in L^2(\mathcal{O})$. 

\item 
The solution of the fully discrete problem \eqref{Frac-discreteh} can be represented by 
\begin{align}\label{unh-rep}
u_n^{(h)}  
&=
 \sum_{j=1}^M \int_{0}^{t_n} E_\tau^{(h)}(t_n-s) \phi_j \bar\partial_\tau W_j(s) \d s  ,
\end{align} 
where the operator $E_\tau^{(h)}(t)$ is given by 
\begin{equation}\label{Etauht-def}
E_\tau^{(h)}(t) \phi:=\!\frac{1}{2\pi {\rm i}}\!\int_{\Gamma_{\theta,\kappa}^{(\tau)}}\!\!\!e^{zt} \frac{z\tau }{e^{z\tau}-1}\delta_\tau(e^{-z\tau})^{\alpha-1}(\delta_\tau(e^{-z\tau})^\alpha 
\! -\! \Delta_h )^{-1}\!P_h\phi\, \d z  
\end{equation}
for $\phi\in L^2(\mathcal{O})$. 
\end{enumerate}

\end{lemma}

The first statement in Lemma \ref{un-rep} has been proved in \cite[Proposition 3.1]{GunzburgerLiWang2017}. The second statement can be proved in the same way, 
replacing the operator $\Delta $ by $\Delta_h$ and $W(\cdot,t)$ by $W^M(\cdot,t)$ (this does not affect the proof therein). 
From Lemma \ref{un-rep}, we see that
\begin{align}\label{error}
u_n-u_n^{(h)}  
=&
\sum_{j=1}^M \int_{0}^{t_n} \big(E_\tau(t_n-s) - E_\tau^{(h)}(t_n-s)\big) \phi_j \bar\partial_\tau W_j(s) \d s
\nonumber \\
&+
\sum_{j=M+1}^\infty \int_{0}^{t_n} E_\tau(t_n-s) \phi_j \bar\partial_\tau W_j(s) \d s 
  \nonumber \\ 
=:& \mathcal{I}_n+\mathcal{J}_n .
\end{align} 
We present the estimates for $\mathcal{I}_n$ and $\mathcal{J}_n$ in subsections \ref{Estimate-In} and \ref{Estimate-Jn}, respectively.

\subsection{Estimate of $\mathcal{I}_n$}\label{Estimate-In}
Now, we start to estimate $\mathcal{I}_n$, i.e., the error of space discretization. The following lemmas are useful in the estimates of $\mathcal{I}_n$ and $\mathcal{J}_n$. 
\begin{lemma}[\cite{Laptev,LY,Strauss2008}]\label{ineq-eigen} 
Let $\mathcal{O}$ denote a bounded domain in $\mathbb{R}^d$, $d\in\{1,2,3\}$. Suppose $\lambda_j$ denotes the $j^{\rm th}$ eigenvalue of the Dirichlet boundary problem for the Laplacian operator $-\Delta$ in $\mathcal{O}$. 
Then, we have 
\begin{align}\label{eigenvalue}
C_0^* j^{\frac{2}{d}} \le \lambda_j\le C_1^* j^{\frac{2}{d}}
\end{align}
for all $j\ge 1$, where the constants $C_0^*$ and $C_1^*$ are independent of $j$. 
\end{lemma}
{\bf Proof.} 
The well-known Weyl's law gives the asymptotic behavior of the eigenvalues of the
Laplacian operator (see \cite{LY} and \cite[pp. 322]{Strauss2008}):
\begin{align}
\lim_{j\rightarrow\infty}\frac{\lambda_j}{j^{2/d}} =(2\pi)^2(B_d |{\cal O}|)^{-\frac{2}{d}}, 
\end{align}
where $B_d$ denotes the volume of the unit $d$-ball. The estimate \eqref{eigenvalue} follows immediately from the above result.\qed

The following lemma is contained in \cite[Lemma 3.2]{GunzburgerLiWang2017}. 
\begin{lemma}\label{ineq-sum}
For any $z\in\Sigma_\varphi$ with $\varphi\in(0,\pi)$, we have
\begin{align}\label{ineq-sum-2}
&\bigg|\frac{1}{z  + \lambda_j }\bigg|
\le \frac{C_\varphi  }{|z|   + \lambda_j } \, ,
\end{align}where $j=1,2,\dots.$
\end{lemma}

The following resolvent estimates will be frequently used in this paper.
\begin{lemma} \label{Lemma:resolvent} 
For $z\in \Sigma_{\theta}$ $($see the definition in Lemma \ref{ineq-1}$)$, with $\theta\in(0,\pi)$, we have the following resolvent estimates:
\begin{align}
&\|(z-\Delta)^{-1}\| + \|(z-\Delta_h)^{-1}\| \le C|z|^{-1},
\label{resolv-est-discr}\\
&\|\Delta^{1-\gamma}(z-\Delta)^{-1}\|+\|\Delta_h^{1-\gamma}(z-\Delta_h)^{-1}\|\le C|z|^{-\gamma},
&&\gamma\in[0,1] . 
\label{resolv-frac-D}
\end{align}
\end{lemma}
{\bf Proof.} 
The first inequality is due to the self-adjointness and nonnegativity of the operators $\Delta$ and $\Delta_h$. These properties guarantee that $\Delta$ and $\Delta_h$ generate bounded analytic semigroup of angle $\frac{\pi}{2}$ on $L^2(\mathcal{O})$ and $(X_h,\|\cdot\|_{L^2(\mathcal{O})})$, respectively; see \cite[Example 3.7.5 and Theorem 3.7.11]{ABHN}.

Recall that $\lambda_j$ and $\phi_j$, $j=1,2,\dots$, are the eigenvalues and eigenfunctions of the operator $-\Delta$ (see Section \ref{sec:notation}). 
The second inequality is due to the interpolation inequality
\begin{align*}
\|\Delta^{1-\gamma}\varphi\| 
&= \Big\|\sum_{j=1}^\infty\lambda_j^{1-\gamma}(\phi_j,\varphi)\phi_j\Big\| \\
&= \Big(\sum_{j=1}^\infty\lambda_j^{2-2\gamma}|(\phi_j,\varphi)|^2\Big)^{\frac{1}{2}}  
 = \Big(\sum_{j=1}^\infty|(\phi_j,\varphi)|^{2\gamma} \lambda_j^{2-2\gamma}|(\phi_j,\varphi)|^{2-2\gamma}\Big)^{\frac{1}{2}} \\
&\le \Big(\sum_{j=1}^\infty|(\phi_j,\varphi)|^2\Big)^{\frac{\gamma}{2}}  
       \Big(\sum_{j=1}^\infty\lambda_j^2|(\phi_j,\varphi)|^2\Big)^{\frac{1-\gamma}{2}} 
       \quad \mbox{(use H\"older's inequality)} \\
&\le \|\varphi\|^{\gamma}\|\Delta\varphi\|^{1-\gamma} .
\end{align*}
Substituting $\varphi=(z-\Delta)^{-1}\phi$ into the inequality above yields 
\begin{align*}
\|\Delta^{1-\gamma}(z-\Delta)^{-1}\phi\| 
&\le \|(z-\Delta)^{-1}\phi\|^{\gamma}\|\Delta(z-\Delta)^{-1}\phi\|^{1-\gamma} \\
&\le (C|z|^{-1}\|\phi\|)^{\gamma}(C\|\phi\|)^{1-\gamma} \\
&\le C|z|^{-\gamma}\|\phi\|.
\end{align*} 
This proves the first part of \eqref{resolv-frac-D}. The estimate of $\|\Delta_h^{1-\gamma}(z-\Delta_h)^{-1}\| $ can be proved similarly (by using eigenvalues and eigenfunctions of $-\Delta_h$).
\qed

The following lemma is concerned with the difference between the continuous and discrete resolvent operators. 

\begin{lemma} \label{ineq-2} 
Let $\alpha\in\big(0,\frac{2}{d}\big)$ and $\delta_\tau(\zeta)$ be defined as in \eqref{definition-d} with the parameters $\kappa$ and $\theta$ satisfying the conditions of Lemma \ref{ineq-1}. Then we have 
\begin{align}\label{ineq-norm}
&\big\|[(\delta_\tau(e^{-z\tau})^\alpha  -\Delta )^{-1}-(\delta_\tau(e^{-z\tau})^\alpha  -\Delta_h )^{-1}P_h]\phi_j\big\| \le Ch^{2\varepsilon}(|z|^{\alpha}+\lambda_j)^{-(1-\varepsilon)}   
\end{align}
for all $\varepsilon\in[0,1]$ and $j=1,2,\dots,M$. 
\end{lemma}
{\bf Proof.} 
First, \eqref{angle-delta} yields that $\delta_\tau(e^{-z\tau})^\alpha\in\Sigma_{\alpha\theta}$ for $z\in\Gamma_{\theta,\kappa}^{(\tau)}$. 
Consequently, we have (cf. \cite[pp. 7]{LubichSloanThomee:1996})
\begin{align}\label{ineq-h2}
\big\|(\delta_\tau(e^{-z\tau})^\alpha  -\Delta )^{-1}-(\delta_\tau(e^{-z\tau})^\alpha  -\Delta_h )^{-1}P_h\big\| \le Ch^2 .
\end{align}

Second, by Lemma \ref{ineq-1} and Lemma \ref{Lemma:resolvent}, there exists a constant $C$ which depends only on $\theta$ and $\alpha$ such that  
\begin{align}\label{ineq-laplacian2}
&
\|(\delta_\tau(e^{-z\tau})^\alpha-\Delta)^{-1}\|
\le
C|\delta_\tau(e^{-z\tau})|^{-\alpha}
\le
C|z|^{-\alpha} &&\forall\, z\in\Gamma_{\theta,\kappa}^{(\tau)},\\
&
\|(\delta_\tau(e^{-z\tau})^\alpha-\Delta_h)^{-1}P_h\|
\le
\|(\delta_\tau(e^{-z\tau})^\alpha-\Delta_h)^{-1}\|\|P_h\|
\le
C|z|^{-\alpha}   &&\forall\, z\in\Gamma_{\theta,\kappa}^{(\tau)} .
\end{align}
Since $(\delta_\tau(e^{-z\tau})^\alpha-\Delta)^{-1}\Delta=(\delta_\tau(e^{-z\tau})^\alpha-\Delta)^{-1}\delta_\tau(e^{-z\tau})^\alpha-I$, it follows that $$\|(\delta_\tau(e^{-z\tau})^\alpha-\Delta)^{-1}\Delta\|\le C.$$ 
Consequently, we have 
\begin{align}\label{ineq-laplacian3}
\|(\delta_\tau(e^{-z\tau})^\alpha-\Delta)^{-1}\phi_j\|
&=
\|(\delta_\tau(e^{-z\tau})^\alpha-\Delta)^{-1}\Delta\Delta^{-1}\phi_j\| \nonumber\\
&\le
\|(\delta_\tau(e^{-z\tau})^\alpha-\Delta)^{-1}\Delta\| \|\Delta^{-1}\phi_j\| \nonumber\\
&\le
C\lambda_j^{-1} ,
\end{align}
which together with $h^2\le C\lambda_M^{-1}$ (since $M=[h^{-d}]+1\sim h^{-d}$ and $\lambda_M\le CM^{\frac{2}{d}}$ by Lemma \ref{ineq-eigen}) implies 
\begin{align}\label{ineq-laplacian4}
\qquad &\|(\delta_\tau(e^{-z\tau})^\alpha-\Delta_h)^{-1} P_h\phi_j \| \nonumber\\
&\le
\|(\delta_\tau(e^{-z\tau}\!)^\alpha\!-\!\Delta)^{-1} \phi_j \| 
+\|[(\delta_\tau(e^{-z\tau}\!)^\alpha\!\!-\!\Delta)^{-1} \!\!-\! (\delta_\tau(e^{-z\tau}\!)^\alpha \!\!-\! \Delta_h)^{-1}P_h]\phi_j\| \nonumber\\
&\le
C\lambda_j^{-1}+Ch^2 
\nonumber\\
&\le
C\lambda_j^{-1}
\end{align}
for $j=1,\cdots,M$. 
Therefore, \eqref{ineq-laplacian2}-\eqref{ineq-laplacian4} leads to
\begin{align}\label{ineq-laplacian5}
\big\|\big[(\delta_\tau(e^{-z\tau})^\alpha  -\Delta )^{-1}-(\delta_\tau(e^{-z\tau})^\alpha  -\Delta_h )^{-1}P_h\big]\phi_j\big\| 
&\le
C\min\{|z|^{-\alpha},\lambda_j^{-1}\} \nonumber\\
&\le 
C(|z|^\alpha+\lambda_j)^{-1}. 
\end{align}

Finally, interpolation between \eqref{ineq-h2} and \eqref{ineq-laplacian5} yields \eqref{ineq-norm}. 
This completes the proof of Lemma \ref{ineq-2}.\qed
\vspace{0.1in}

Now, we turn to the estimate of $\mathcal{I}_n$. From Lemma \ref{ineq-2}, and choosing $\beta\in(0,1)$, it is easy to derive 
\begin{align*}
&\|\big(E_\tau(t) - E_\tau^{(h)}(t)\big) \phi_j\|^2 \\
&\le
\bigg(\int_{\Gamma_{\theta,\kappa}^{(\tau)}} |e^{zt}| |\delta_\tau(e^{-z\tau})|^{\alpha-1}\|[(\delta_\tau(e^{-z\tau})^{\alpha}-\Delta)^{-1}-(\delta_\tau(e^{-z\tau})^\alpha-\Delta_h)^{-1}P_h ]\phi_j \|\d z \bigg)^2 \\
&\le 
\bigg(C\int_{\Gamma_{\theta,\kappa}^{(\tau)}}e^{|z|t\cos({\rm arg}(z))} |z|^{\alpha-1}
(|z|^{\alpha}+\lambda_j)^{-(1-\varepsilon)}  h^{2\varepsilon} \, |\d z| \bigg)^2 \\
&\le
Ch^{4\varepsilon}\bigg(\int_{\Gamma_{\theta,\kappa}^{(\tau)}} e^{|z|t\cos(\arg(z))}\frac{|\d z|}{|z|^\beta}\bigg)
\left(\int_{\Gamma_{\theta,\kappa}^{(\tau)}} e^{|z|t\cos(\arg(z))} \frac{|z|^{2\alpha-2+\beta}}{(|z|^\alpha+\lambda_j)^{2-2\varepsilon}}|\d z| \right) \\
&\le
Ch^{4\varepsilon} t^{\beta-1}\int_{\Gamma_{\theta,\kappa}^{(\tau)}} e^{|z|t\cos(\arg(z))}\frac{|z|^{2\alpha-2+\beta}}{(|z|^\alpha+\lambda_j)^{2-2\varepsilon}}|\d z| \, ,
\end{align*}
where we use $\kappa= \frac{1}{T}$ and 
\begin{align}\label{integral-gamma}
\int_{\Gamma_{\theta,\kappa}^{(\tau)}} e^{|z|t\cos(\arg(z))}\frac{|\d z|}{|z|^\beta}
&\le
\int_\kappa^{\frac{\pi}{\tau\sin(\theta)}}e^{-rt|\cos(\theta)|}\frac{\d r}{r^\beta}
+\int_{-\theta}^\theta e^{\kappa t\cos(\varphi)} \frac{\d\varphi}{\kappa^{\beta-1}} \nonumber\\
&\le
Ct^{\beta-1}. 
\end{align}
Since $\frac{W_j(t_i)-W_j(t_{i-1})}{\tau}$, $i=1,\dots,n$, and $j=1,2,\dots$, are stochastically independent of each other, it follows that 
\begin{align}\label{ineq-In}
{\mathbb E}\|\mathcal{I}_n\|^2
&= 
\sum_{j=1}^M {\mathbb E} \bigg\|\sum_{i=1}^n\frac{W_j(t_i)-W_j(t_{i-1})}{\tau} \int_{t_{i-1}}^{t_i}(E_\tau(t_n-s)-E_\tau^{(h)}(t_n-s))\phi_j \d s  \bigg\|^2 \nonumber \\
&=  \sum_{j=1}^M  \sum_{i=1}^n\frac{1}{\tau} \bigg\|\int_{t_{i-1}}^{t_i}(E_\tau(t_n-s)-E_\tau^{(h)}(t_n-s))\phi_j \d s  \bigg\|^2 \nonumber \\ 
&\le  \sum_{j=1}^M  \sum_{i=1}^n  \int_{t_{i-1}}^{t_i}\| E_\tau(t_n-s)-E_\tau^{(h)}(t_n-s))\phi_j\|^2 \d s   \nonumber 
\quad\mbox{(H\"older's inequality)}\\ 
&=
\sum_{j=1}^M \int_{0}^{t_n} \|\big(E_\tau(t) - E_\tau^{(h)}(t)\big) \phi_j \|^2 
\d t  \nonumber \\ 
&\le 
Ch^{4\varepsilon}\int_{0}^{t_n}t^{\beta-1}   \int_{\Gamma_{\theta,\kappa}^{(\tau)}} 
e^{|z|t\cos({\rm arg}(z))} \sum_{j=1}^M\frac{|z|^{2\alpha-2+\beta}}{(|z|^{\alpha}+\lambda_j)^{2-2\varepsilon}}   \, |\d z|  \d t \nonumber \\
&= 
Ch^{4\varepsilon}\int_{0}^{t_n}  t^{\beta-1}   \int_{\Gamma_{\theta,\kappa}^{(\tau)}} 
e^{|z|t\cos({\rm arg}(z))} |z|^{2\varepsilon\alpha-2+\beta} \sum_{j=1}^M\bigg(\frac{|z|^{\alpha}}{|z|^{\alpha}+\lambda_j}\bigg)^{2-2\varepsilon} \, |\d z|  \d t  .
\end{align}

Since $\big(\frac{|z|^{\alpha}}{|z|^{\alpha}+\lambda_j}\big)^{2-2\varepsilon}\le 1$ for $\varepsilon\in[0,1]$, it follows that 
$$
\sum_{j=1}^M\bigg(\frac{|z|^{\alpha}}{|z|^{\alpha}+\lambda_j}\bigg)^{2-2\varepsilon}\le M\le Ch^{-d} ,
$$
where the last inequality is due to our choice $M=[h^{-d}]+1$. 
Thus \eqref{ineq-In} reduces to  
\begin{align} 
\mathbb{E}\|\mathcal{I}_n\|^2
&\le 
Ch^{4\varepsilon-d} \int_{0}^{t_n}  t^{\beta-1}    \int_{\Gamma_{\theta,\kappa}^{(\tau)}} 
e^{|z|t\cos({\rm arg}(z))} |z|^{2\varepsilon\alpha-2+\beta} \, |\d z|  \d t \nonumber \\
&\le 
Ch^{4\varepsilon-d} \int_{0}^{t_n} t^{\beta-1}   
\int_{\kappa}^{\frac{\pi}{\tau\sin(\theta)}} e^{-rt|\cos(\theta)|} r^{2\varepsilon\alpha-2+\beta} \, \d r \d t \nonumber\\
&\quad
+
Ch^{4\varepsilon-d} \int_0^{t_n} t^{\beta-1}\int_{-\theta}^{\theta} e^{\kappa t\cos(\varphi)}   \kappa^{2\varepsilon\alpha-1+\beta} \, \d\varphi\d t  \nonumber\\ 
&= 
Ch^{4\varepsilon-d} 
\int_{\kappa}^{\frac{\pi}{\tau\sin(\theta)}}r^{2\varepsilon\alpha-2+\beta}  \int_{0}^{t_n} t^{\beta-1}    e^{-rt|\cos(\theta)|} \,\d t \d r \nonumber\\
&\quad
+Ch^{4\varepsilon-d}  
\int_{0}^{t_n} t^{\beta-1}\kappa^{2\varepsilon\alpha-1+\beta}    \int_{-\theta}^{\theta} e^{\kappa t\cos(\varphi)}   \, \d\varphi \d t \nonumber\\
&\le 
Ch^{4\varepsilon-d} 
\bigg(\int_{\kappa}^{\frac{\pi}{\tau\sin(\theta)}}  r^{2\varepsilon\alpha-2} \, \d r
+\kappa^{2\varepsilon\alpha-1+\beta}\int_{0}^{t_n} t^{\beta-1}   \d t\bigg)   \nonumber\\
&\le 
C h^{4\varepsilon-d} \bigg(\frac{1}{1-2\varepsilon\alpha}\kappa^{2\varepsilon\alpha-1}-\frac{1}{1-2\varepsilon\alpha}\bigg(\frac{\pi}{\tau\sin(\theta)}\bigg)^{2\varepsilon\alpha-1}
+ 
\kappa^{2\varepsilon\alpha-1} (\kappa t_n)^\beta\bigg)  
. \nonumber
\end{align} 
For $\alpha\in[\frac12,\frac{d}{2})$, we choose $\varepsilon=\frac{1-1/\ell_h}{2\alpha}$ with $\ell_h=\ln(e+1/h)$ and $h^{-1/\ell_h}\le e$ (can recall that $\kappa= \frac{1}{T}$),
we have  
\begin{align}\label{In--1}
\mathbb{E}\|\mathcal{I}_n\|^2
&\le 
CT^{\frac{1}{\ell_h}} \ell_h h^{\frac{2}{\alpha}-d},\quad\mbox{for}\,\,\, \alpha\in\Big[\frac12,\frac{d}{2}\Big) .
\end{align} 
For $\alpha\in(0,\frac{1}{2})$, we choose $\varepsilon=1$ and get 
\begin{align}\label{In--2}
\mathbb{E}\|\mathcal{I}_n\|^2
&\le
CT^{1-2\alpha}  h^{4-d} ,\quad\mbox{for}\,\,\, \alpha\in\Big(0,\frac12\Big) .
\end{align}

\subsection{Estimate of $\mathcal{J}_n$}\label{Estimate-Jn}
In this subsection, we present the estimate of $\mathcal{J}_n$ in \eqref{error}, completing the proof of \eqref{error-u-unh}. 
In view of the definition of $E_\tau(t)$ in \eqref{eqn:EF2}, by using Lemma \ref{ineq-1} and \eqref{ineq-sum-2} to estimate $|\delta_\tau(e^{-z\tau})^{\alpha-1}|$ and $|(\delta_\tau(e^{-z\tau})^\alpha  + \lambda_j )^{-1}|$, we obtain
\begin{align} \label{E_tautphij}
\| E_\tau(t)\phi_j\|^2 
&=
\bigg|\frac{1}{2\pi i}\int_{\Gamma_{\theta,\kappa}^{(\tau)}}e^{zt} \frac{z\tau }{e^{z\tau}-1}\delta_\tau(e^{-z\tau})^{\alpha-1}(\delta_\tau(e^{-z\tau})^\alpha  + \lambda_j )^{-1} \d z \bigg|^2  \nonumber\\
&\le
C \bigg( \int_{\Gamma_{\theta,\kappa}^{(\tau)}} |e^{zt}| \bigg|\frac{z\tau }{e^{z\tau}-1}\bigg| |\delta_\tau(e^{-z\tau})|^{\alpha-1}(|\delta_\tau(e^{-z\tau})|^\alpha  + \lambda_j )^{-1} \d z \bigg)^2   \nonumber\\
&\le
C \bigg( \int_{\Gamma_{\theta,\kappa}^{(\tau)}} |e^{zt}|   |z|^{\alpha-1}(|z|^\alpha  + \lambda_j )^{-1} \d z \bigg)^2  \quad\mbox{(use \eqref{1-expztau-1} here)} \nonumber\\
&\le 
C \bigg(\int_{\Gamma_{\theta,\kappa}^{(\tau)}}e^{|z|t\cos({\rm arg}(z))} \frac{|\d z|}{|z|^{\beta}} |\d z| \bigg)
\bigg(\int_{\Gamma_{\theta,\kappa}^{(\tau)}}e^{|z|t\cos({\rm arg}(z))} \frac{|z|^{2\alpha-2+\beta}}{(|z|^\alpha  + \lambda_j)^2} |\d z| \bigg) \nonumber\\
&\le 
Ct^{\beta-1} \int_{\Gamma_{\theta,\kappa}^{(\tau)}}e^{|z|t\cos({\rm arg}(z))} \frac{|z|^{2\alpha-2+\beta}}{(|z|^{\frac{d\alpha}{2}}  + j)^{\frac{4}{d}}} |\d z| ,
\end{align} 
where we have used \eqref{integral-gamma}. In view of \eqref{error}, we have 
\begin{align} \label{ineq-Jn}
{\mathbb E}\|\mathcal{J}_n\|^2
&\le 
\int_0^{t_n}  \sum_{j=M+1}^\infty \| E_\tau(t)\phi_j\|^2 \d t  \nonumber\\
&\le 
C\int_0^{t_n}   t^{\beta-1} \int_{\Gamma_{\theta,\kappa}^{(\tau)}}e^{|z|t\cos({\rm arg}(z))}\sum_{j=M+1}^\infty\frac{|z|^{2\alpha-2+\beta}}{(|z|^{\frac{d\alpha}{2}}  + j)^{\frac{4}{d}}} |\d z|\d t \nonumber\\
&\le 
C\int_0^{t_n}   t^{\beta-1} \int_{\Gamma_{\theta,\kappa}^{(\tau)}}e^{|z|t\cos({\rm arg}(z))}\frac{|z|^{2\alpha-2+\beta}}{(|z|^{\frac{d\alpha}{2}}  +M)^{\frac{4}{d}-1}} |\d z|\d t \nonumber\\
&\le 
C\int_0^{t_n}t^{\beta-1}  \int_{\kappa}^{\frac{\pi}{\tau\sin(\theta)}}e^{-rt|\cos(\theta)|}\frac{r^{2\alpha-2+\beta}}{(r^{\frac{d\alpha}{2}}  +M)^{\frac{4}{d}-1}} \d r\d t \nonumber\\
&\quad
+C\int_0^{t_n}   t^{\beta-1}\int_{-\theta}^{\theta}  e^{\kappa t\cos(\varphi)}\frac{\kappa^{2\alpha-1+\beta}}{(\kappa^{\frac{d\alpha}{2}}  +M)^{\frac{4}{d}-1}}\d \varphi \d t  \nonumber\\
&=
C\int_{\kappa}^{\frac{\pi}{\tau\sin(\theta)}} \frac{r^{2\alpha-2+\beta}}{(r^{\frac{d\alpha}{2}}  +M)^{\frac{4}{d}-1}}\int_0^{t_n}t^{\beta-1}e^{-rt|\cos(\theta)|} \d t\d r \nonumber\\
&\quad
+C\int_0^{t_n}   t^{\beta-1}\frac{\kappa^{2\alpha-1+\beta}}{(\kappa^{\frac{d\alpha}{2}}  +M)^{\frac{4}{d}-1}}\int_{-\theta}^{\theta}  e^{\kappa t\cos(\varphi)}\d \varphi \d t  \nonumber\\
&\le 
C \int_{\kappa}^{\frac{\pi}{\tau\sin(\theta)}}\frac{r^{2\alpha-2}}{(r^{\frac{d\alpha}{2}}  +M)^{\frac{4}{d}-1}} \d r
+ \frac{C(\kappa t_n)^\beta\kappa^{2\alpha-1}}{(\kappa^{\frac{d\alpha}{2}}  +M)^{\frac{4}{d}-1}} \nonumber \\
&\le 
C \bigg(\int_{\kappa}^{\min\{M^{\frac{2}{d\alpha}} ,\frac{\pi}{\tau\sin(\theta)}\}}\frac{r^{2\alpha-2}}{M^{\frac{4}{d}-1}} \d r
+\int_{\min\{M^{\frac{2}{d\alpha}} ,\frac{\pi}{\tau\sin(\theta)} \} }^{\frac{\pi}{\tau\sin(\theta)}}\frac{r^{2\alpha-2}}{r^{2\alpha-\frac{d\alpha}{2}}} \d r\bigg) \nonumber\\
&\quad
+\frac{C  \kappa^{2\alpha-1}}{(\kappa^{\frac{d\alpha}{2}}  +M)^{\frac{4}{d}-1}} \nonumber\\
&\le 
C 
\bigg(\int_{\kappa}^{\min\{M^{\frac{2}{d\alpha}} ,\frac{\pi}{\tau\sin(\theta)} \} }\frac{r^{2\alpha-2}}{M^{\frac{4}{d}-1}} \d r
+M^{1-\frac{2}{d\alpha}} \bigg)
+ \frac{C \kappa^{2\alpha-1}}{(\kappa^{\frac{d\alpha}{2}}  +M)^{\frac{4}{d}-1}}  .
\end{align} 
We estimate $\mathbb{E}\|\mathcal{J}_n\|^2$ in three different cases by using 
$M=[h^{-d}]+1\le Ch^{-d}$. 

In the case $\alpha\in(\frac{1}{2},\frac{2}{d})$, 
\eqref{ineq-Jn} reduces to
\begin{align}\label{Jn--1}
{\mathbb E}\|\mathcal{J}_n\|^2 
&\le 
C 
\bigg(\int_{\kappa}^{M^{\frac{2}{d\alpha}} }\frac{r^{2\alpha-2}}{M^{\frac{4}{d}-1}} \d r
+M^{1-\frac{2}{d\alpha}} \bigg) 
+ CM^{1-\frac{4}{d}} \nonumber \\
&\le
CM^{1-\frac{2}{d\alpha}}\le Ch^{\frac{2}{\alpha}-d}  .
\end{align}

In the case $\alpha\in(0,\frac{1}{2})$, \eqref{ineq-Jn} yields 
\begin{align}\label{Jn--2}
{\mathbb E}\|\mathcal{J}_n\|^2 
\le C \big(T^{1-2\alpha}M^{1-\frac{4}{d}}
+M^{1-\frac{2}{d\alpha}}  \big) 
&\le 
C(1+T^{1-2\alpha}) M^{1-\frac{4}{d}} \nonumber\\
&\le C(1+T^{1-2\alpha}) h^{4-d}   . 
\end{align}

In the case $\alpha=\frac{1}{2}$, \eqref{ineq-Jn} implies  
\begin{align}\label{Jn--3}
{\mathbb E}\|\mathcal{J}_n\|^2 
\le
C 
\bigg(\int_{\kappa}^{M^{\frac{2}{d\alpha}}}\frac{r^{-1}}{M^{\frac{4}{d}-1}} \d r
+M^{1-\frac{2}{d\alpha}} \bigg)
+ CM^{1-\frac{4}{d}}  
&\le 
C (\ln M+\ln T)M^{1-\frac{4}{d}} \nonumber\\
&\le C(\ell_h+\ln T) h^{\frac{2}{\alpha}-d} .
\end{align}

It is well-known that $(\mathbb{E}\|u_n-u_n^{(h)}\|)^2\le\mathbb{E}\|u_n-u_n^{(h)}\|^2$, 
substituting \eqref{In--1}-\eqref{In--2} and \eqref{Jn--1}-\eqref{Jn--3} into \eqref{error} yields \eqref{space-discr}, completing the proof of \eqref{error-u-unh}.

\section{Deterministic problem: 
estimate of $\|v(\cdot,t_n)-v_n^{(h)}\|$ 
}\label{sec:deter}

In this section, we estimate the error $\|v(\cdot,t_n)-v_n^{(h)}\|$ by minimizing the regularity requirement on $\psi_0$ and $f$ to match the convergence rate proved in the last section, where $v$ and $v_n^{(h)}$ are the solutions of \eqref{Deter-SPDE2}  and \eqref{Deter-discreteh} or \eqref{Deter-discreteh-2}, respectively. In particular, we prove the following estimate for $\alpha\in\big(0,\frac{2}{d}\big)$:
\begin{align}\label{error-v-vnh}
\|v(\cdot,t_n)-v_n^{(h)}\|\le 
C\big(\|\psi_0\|_{\dot H^{\chi}(\mathcal{O})} \!+\! \|f\|_{L^{\frac{4}{2+\alpha d},1}(0,t_n;L^2(\mathcal{O}))}\big) 
\big(\tau^{\frac{1}{2}-\frac{\alpha d}{4}} +h^{\chi}\big) ,
\end{align}
where $\chi\in(0,2-\frac{d}{2})$ is defined in \eqref{Def-chi}. 
To this end, we introduce the semi-discrete finite element problem
\begin{equation}\label{semi-FEM}
\left\{\begin{aligned}
&\partial_t v^{(h)}  -\Delta_h \partial_t^{1-\alpha} v^{(h)}
= P_hf  \\
&v^{(h)}(\cdot,0)=P_h\psi_0 .
\end{aligned}
\right.
\end{equation}
Then \eqref{Deter-discreteh} or \eqref{Deter-discreteh-2}  can be viewed as the time discretization of \eqref{semi-FEM},  
with the following decomposition:
\begin{align}\label{error-split-FEM}
\|v(\cdot,t_n)-v_n^{(h)}\|
&\le
\|v(\cdot,t_n)-v^{(h)}(\cdot,t_n)\| + \|v^{(h)}(\cdot,t_n)-v_n^{(h)}\| . 
\end{align}
We estimate $\|v(\cdot,t_n)-v^{(h)}(\cdot,t_n)\|$ and $\|v^{(h)}(\cdot,t_n)-v_n^{(h)}\|$ in the next two subsections, respectively.

\subsection{Spatial discretization: estimate of $\|v(\cdot,t_n)-v^{(h)}(\cdot,t_n)\|$} 
\label{sd-f} 
In this subsection we estimate $\|v(\cdot,t_n)-v^{(h)}(\cdot,t_n)\|$, where 
$v$ and $v^{(h)}$ are solutions of \eqref{Deter-SPDE2} and \eqref{semi-FEM}, respectively. 
Next, we consider the following three cases:\\
Case 1: $\psi_0\neq 0$,  $f=0$, $\alpha\in(0,1]\cap(0,\frac{2}{d})$.\\ 
Case 2: $\psi_0\neq 0$,  $f=0$, $\alpha\in(1,\frac{2}{d})$. 
(This requires $d=1$, cf. Theorem \ref{MainTHM})\\
Case 3: $\psi_0= 0$, $f\neq 0$. 

In Case 1, the solutions of \eqref{Deter-SPDE2} and \eqref{semi-FEM} are given by 
$v(\cdot,t_n)=\psi_0$ and $v^{(h)}(\cdot,t_n)=P_h\psi_0$, $n=1,2,\dots$, respectively 
(this follows from setting $f=0$ and $\sigma=0$ in \eqref{Mild-sol-} of Appendix \ref{AppendixA}). 
Consequently, we have
\begin{align}
\|v(\cdot,t_n) - v^{(h)} (\cdot,t_n)  \|
=\|\psi_0-P_h\psi_0\| 
\le C \|\psi_0\|_{\dot H^{\chi}(\mathcal{O})} h^{\chi} ,
\end{align}
where we have used \eqref{StabPh-} in the last inequality. 

In Case 2, we note that the PDE problem \eqref{Deter-SPDE2} is equivalent to (multiplying both sides by the operator $\partial_t^{\alpha-1}$)  
\begin{align}\label{Deter-SPDE2-2}
\left\{
\begin{array}{ll}
\partial_t^\alpha  v(x,t) 
-\Delta v(x,t) = 0
&\quad  (x,t)\in \mathcal{O}\times {\mathbb R}_+  \\
v(x,t)=0   &\quad (x,t)\in \partial\mathcal{O}\times {\mathbb R}_+    \\
v(x,0)=\psi_0(x) \,\,\,\mbox{and}\,\,\, \partial_t v(x,0)=0    &\quad x\in \mathcal{O}   
\end{array}\right.
\end{align}
and the finite element problem \eqref{semi-FEM} is equivalent to (for the same reason) 
\begin{equation}\label{Deter-FEM2-2}
\left\{\begin{aligned}
&\partial_t^\alpha v^{(h)} 
-\Delta_h v^{(h)}
= 0 \\
&v^{(h)}(\cdot,0)=P_h\psi_0 
\,\,\,\mbox{and}\,\,\,
\partial_tv^{(h)}(\cdot,0)=0 .
\end{aligned}
\right.
\end{equation}
By Laplace transform, $v=E(t)\psi_0$ and $v^{(h)}=E^{(h)}(t)P_h\psi_0$, where the operators $E(t):L^2({\cal O})\rightarrow L^2({\cal O})$ and $E^{(h)}(t):X_h\rightarrow X_h$ are given in Appendix \ref{AppendixA}. 
An error estimate for \eqref{Deter-SPDE2-2} and \eqref{Deter-FEM2-2} was presented in  
\cite[Theorems 3.2]{JinLazarovZhouSISC2016}:   
\begin{align}\label{space-vh1}
\|v(\cdot,t_n)-v^{(h)}(\cdot,t_n)\|  \le C \|\psi_0\|_{\dot H^{2}(\mathcal{O})} h^2.
\end{align}
The boundedness of the solution operators $E(t):L^2(\mathcal{O})\rightarrow L^2(\mathcal{O})$ and $E^{(h)}(t):X_h\rightarrow X_h$ (see Appendix \ref{AppendixA}, \eqref{BD-Et} and \eqref{BD-Eht}) implies
\begin{align}\label{space-vh2}
\begin{aligned}
\|v(\cdot,t_n)-v^{(h)}(\cdot,t_n) \|
&=\|E(t_n)\psi_0-E^{(h)}(t_n)P_h\psi_0 \| \\
&\le \|E(t_n)\psi_0\|+\|E^{(h)}(t_n)P_h\psi_0 \| \\
&\le C\|\psi_0\|  .
\end{aligned}
\end{align}
Then the interpolation between \eqref{space-vh1} and \eqref{space-vh2} yields 
\begin{align}\label{space-vhf}
\begin{aligned}
\|v(\cdot,t_n)-v^{(h)}(\cdot,t_n) \|
&\le C\|\psi_0\|_{\dot H^{\chi}(\mathcal{O})} h^{\chi}   .
\end{aligned}
\end{align}

In Case 3, we have 
\begin{align}\label{Phv-vht}
P_h  v(\cdot,t_n) - v^{(h)} (\cdot,t_n)  
&=P_h  \bigg(v(\cdot,t_n) - v^{(h)} (\cdot,t_n)\bigg)  \nonumber \\
&=
P_h\bigg(\int_{0}^{t_n} E(t_n-s) f (\cdot,s) \d s - \int_{0}^{t_n}E^{(h)}(t_n-s)P_hf (\cdot,s)\d s\bigg) \nonumber\\
&=
P_h \int_{0}^{t_n} E(t_n-s) (f (\cdot,s)-P_hf (\cdot,s)) \d s \nonumber\\
&\quad + P_h \int_{0}^{t_n}\Big(E(t_n-s) - E^{(h)}(t_n-s)P_h\Big) P_hf (\cdot,s)\d s \nonumber\\
&=: \mathcal{I}^{(h)}(t_n)+\mathcal{J}^{(h)}(t_n) ,
\end{align} 
where
\begin{align} \label{space-Iht}
\begin{aligned}
\mathcal{I}^{(h)}(t_n)
&=\int_0^{t_n}\frac{1}{2\pi\mathrm{i}}\int_{\Gamma_{\theta,\kappa}} 
     e^{z(t_n-s)}z^{\alpha-1}P_h  (z^\alpha - \Delta)^{-1}  (f (\cdot,s) -P_hf (\cdot,s))\d z\d s \\
&=\int_0^{t_n}\frac{1}{2\pi\mathrm{i}}\int_{\Gamma_{\theta,\kappa}} 
     e^{z(t_n-s)}z^{\alpha-1} P_h w_z(\cdot,s)\d z \d s 
     \end{aligned}
\end{align} 
with 
\begin{align}\label{Def-w-int-f}
w_z(\cdot,s)=(z^\alpha - \Delta)^{-1}(f(\cdot,s)-P_hf(\cdot,s)) ,
\end{align}  
which satisfies 
\begin{align} \label{H2w-it-f}
\|w_z(\cdot,s)\|_{\dot H^2(\mathcal{O})}
&=    \|\Delta^{-1}\Delta(z^\alpha - \Delta)^{-1}(f(\cdot,s)-P_hf(\cdot,s))\|_{\dot H^2(\mathcal{O})} \nonumber \\
&\le C \| \Delta(z^\alpha - \Delta)^{-1}(f(\cdot,s)-P_hf(\cdot,s))\| \nonumber \\
&\le C \|f(\cdot,s)-P_hf(\cdot,s)\|.
\end{align} 

By using the Ritz projection operator $R_h$, equation \eqref{Def-w-int-f} implies 
\begin{align*} 
((z^\alpha - \Delta_h)P_hw_z+\Delta_h(P_hw_z-R_hw_z),\phi_h)  
&=z^\alpha (w_z, \phi_h)+(\nabla  R_h w_z , \nabla \phi_h) \\
&=z^\alpha (w_z, \phi_h)+(\nabla w_z , \nabla \phi_h)  \\
&=((z^\alpha - \Delta)w_z,\phi_h) \\
&=(f-P_hf,\phi_h)=0  \quad\forall\, \phi_h\in X_h ,
\end{align*} 
i.e., $P_hw_z=-(z^\alpha - \Delta_h)^{-1}\Delta_h(P_hw_z-R_hw_z)$. Consequently, we obtain
\begin{align*} 
|z|^{\gamma\alpha}\|P_hw_z(\cdot,s)\| 
&= \||z|^{\gamma\alpha}(z^\alpha - \Delta_h)^{-\gamma}[\Delta_h(z^\alpha - \Delta_h)^{-1}]^{1-\gamma}\Delta_h^\gamma (P_hw_z(\cdot,s)-R_hw_z(\cdot,s))\| \\
&\le C\|\Delta_h^\gamma (P_hw_z(\cdot,s)-R_hw_z(\cdot,s))\| \\
&\le C\|w_z(\cdot,s)\|_{\dot H^2(\mathcal{O})} h^{2-2\gamma} \qquad \qquad\mbox{(here we use \eqref{StabPh4})}\\
&\le C\|f(\cdot,s)-P_hf(\cdot,s)\| h^{2-2\gamma} 
 \qquad \mbox{(here we use \eqref{H2w-it-f})}\\
&\le  C\|f(\cdot,s)\| h^{2-2\gamma}.
\end{align*}
By choosing $\gamma=1-\frac{1}{2}\chi$ so that 
$2-2\gamma=\chi$, 
the inequality \eqref{space-Iht} reduces to 
\begin{align}\label{space-Ihf}
\begin{aligned}
\|\mathcal{I}^{(h)}(t_n)\|
&\le C\int_0^{t_n}\int_{\Gamma_{\theta,\kappa}} 
          |e^{z(t_n-s)}||z|^{(1-\gamma)\alpha-1} |z|^{\gamma\alpha}\|P_h w(\cdot,s)\| |\d z| \d s  \\
&\le C\int_0^{t_n} \int_{\Gamma_{\theta,\kappa}}  
          |e^{z(t_n-s)}| |z|^{\frac{1}{2}\chi\alpha -1}\|f(\cdot,s)\|h^{\chi} |\d z| \d s\\
&\le Ch^{\chi} \int_0^{t_n} (t_n-s)^{-\frac{1}{2}\chi\alpha} \|f(\cdot,s)\|  \d s .
\end{aligned}
\end{align} 
Furthermore, applying the similar analysis in Lemma \ref{ineq-2}, we have 
\begin{align*}
&\|(z^\alpha - \Delta)^{-1}-(z^\alpha - \Delta_h)^{-1}P_h\|\le  Ch^2  ,
&&\mbox{(cf. \cite[pp. 7]{LubichSloanThomee:1996})} \\
&\|(z^\alpha - \Delta)^{-1}-(z^\alpha - \Delta_h)^{-1}P_h\|\le  C|z|^{-\alpha} .
&&\mbox{(resolvent estimate, see \eqref{resolv-est-discr})}
\end{align*}
The interpolation of the last two inequalities yields 
\begin{align*}
&\|(z^\alpha - \Delta)^{-1}-(z^\alpha - \Delta_h)^{-1}P_h\| 
\le  C|z|^{-\gamma\alpha} h^{2-2\gamma} \quad\forall\, \gamma\in[0,1].
\end{align*}
Again, by choosing $\gamma=1-\frac{1}{2}\chi$ (so that $2-2\gamma=\chi$), we have 
\begin{align} \label{space-Jhf}
&\|\mathcal{J}^{(h)}(t_n) \|  \nonumber\\
& = \bigg\|P_h\int_0^{t_n}\frac{1}{2\pi\mathrm{i}}\int_{\Gamma_{\theta,\kappa}} 
     \!\! e^{z(t_n-s)}z^{\alpha-1}  \Big((z^\alpha - \Delta)^{-1} \!\! - \! (z^\alpha - \Delta_h)^{-1}P_h\Big)P_hf (\cdot,s)\d z\d s\bigg\| \nonumber \\
& \le \int_0^{t_n}\int_{\Gamma_{\theta,\kappa}} 
     |e^{z(t_n-s)}| |z|^{\alpha-1}  \|(z^\alpha - \Delta)^{-1}-(z^\alpha - \Delta_h)^{-1}P_h\| 
     \|P_hf (\cdot,s)\| |\d z| \d s \nonumber \\
& \le C\int_0^{t_n}\int_{\Gamma_{\theta,\kappa}} 
     |e^{z(t_n-s)}| |z|^{\alpha-1}  |z|^{-(1-\frac{1}{2}\chi)\alpha} h^{\chi}
     \|f (\cdot,s)\| |\d z| \d s \nonumber \\
& \le Ch^{\chi} \int_0^{t_n} \bigg(\int_{\Gamma_{\theta,\kappa}} 
     |e^{z(t_n-s)}| |z|^{\frac{1}{2}\chi\alpha-1} 
      |\d z| \bigg)\|f (\cdot,s)\|\d s \nonumber \\
& \le Ch^{\chi} \int_0^{t_n} (t_n-s)^{- \frac{1}{2}\chi\alpha} \|f (\cdot,s)\|\d s .
\end{align} 
Substituting \eqref{space-Ihf} and \eqref{space-Jhf} into \eqref{Phv-vht} yields  
\begin{align}\label{space-Phv-vh}
\|P_h  v(\cdot,t_n) - v^{(h)} (\cdot,t_n)  \|
\le Ch^{\chi} \int_0^{t_n} (t_n-s)^{- \frac{1}{2}\chi\alpha} \|f (\cdot,s)\|\d s .
\end{align} 
Furthermore, we can see that
\begin{align}\label{space-v-Phv}
&\|v(\cdot,t_n)-P_h  v(\cdot,t_n) \| \nonumber\\
&\le Ch^{\chi} \|v(\cdot,t_n)\|_{\dot H^{\chi} (\mathcal{O})}  \nonumber\\
&\le Ch^{\chi} \|\Delta^{\frac{1}{2}\chi}v(\cdot,t_n)\|  \nonumber \\
&= 
Ch^{\chi} \bigg\|\int_0^{t_n} \frac{1}{2\pi\mathrm{i}}\int_{\Gamma_{\theta,\kappa}} 
e^{z(t_n-s)}z^{\alpha-1}\Delta^{\frac{1}{2}\chi}(z^\alpha-\Delta)^{-1} f(\cdot,s)\d z \d s \bigg\|  \nonumber\\
&\le 
Ch^{\chi} \int_0^{t_n} \int_{\Gamma_{\theta,\kappa}} 
|e^{z(t_n-s)}| |z|^{\alpha-1}
\|\Delta^{\frac{1}{2}\chi}(z^\alpha-\Delta)^{-1}\| 
\|f(\cdot,s)\| |\d z| \d s  \nonumber\\
&\le 
Ch^{\chi} \int_0^{t_n} \int_{\Gamma_{\theta,\kappa}} 
|e^{z(t_n-s)}| |z|^{\frac{1}{2}\chi\alpha-1} 
\|f(\cdot,s)\| |\d z| \d s \quad(\textrm{here we use }\eqref{resolv-frac-D}) \nonumber\\
&\le Ch^{\chi} \int_0^{t_n} (t_n-s)^{-\frac{1}{2}\chi\alpha} \|f(\cdot,s)\|  \d s . 
\end{align} 
The estimates \eqref{space-Phv-vh}-\eqref{space-v-Phv} imply 
\begin{align}
\|v(\cdot,t_n) - v^{(h)} (\cdot,t_n)  \|
&\le Ch^{\chi} \int_0^{t_n} (t_n-s)^{-\frac{1}{2}\chi\alpha} \|f (\cdot,s)\|\d s  \nonumber \\
&\le Ch^{\chi}\|f \|_{L^{\frac{2}{2-\chi\alpha},1}(0,T;L^2(\mathcal{O}))}, 
\quad\mbox{(here we use \eqref{Lp1-use})}
\end{align} 
completing the proof in Case 3. 

The combination of Cases 1, 2 and 3 yields 
\begin{align}\label{vvvvvv}
\|v(\cdot,t_n) - v^{(h)} (\cdot,t_n)  \|
&\le C 
\big(\|\psi_0\|_{\dot H^{\chi}(\mathcal{O})} 
+\|f \|_{L^{\frac{2}{2-\chi\alpha},1}(0,T;L^2(\mathcal{O}))}\big) h^{\chi} \nonumber \\
&\le C 
\big(\|\psi_0\|_{\dot H^{\chi}(\mathcal{O})} 
+\|f \|_{L^{\frac{4}{2+\alpha d},1}(0,T;L^2(\mathcal{O}))}\big) h^{\chi}  ,
\end{align}
where we have used the fact $\frac{2}{2-\chi\alpha}\le \frac{4}{2+\alpha d}$ in the last inequality.

\subsection{Temporal discretization: estimate of $\|v^{(h)}(\cdot,t_n)-v_n^{(h)}\|$}
\label{td-f} 
To estimate $\|v^{(h)}(\cdot,t_n)-v_n^{(h)}\|$, we consider the following three cases:\\
Case 1': $\psi_0\neq 0$,  $f=0$, $\alpha\in(0,1]\cap(0,\frac{2}{d})$.\\ 
Case 2': $\psi_0\neq 0$,  $f=0$, $\alpha\in(1,\frac{2}{d})$. 
(This requires $d=1$, cf. Theorem \ref{MainTHM})\\
Case 3': $\psi_0= 0$, $f\neq 0$.  

In Case 1', it is straightforward to verify that the solutions of \eqref{semi-FEM} and \eqref{Deter-discreteh} 
are given by 
$v^{(h)}(\cdot,t_n)=v^{(h)}_n=P_h\psi_0$, $n=1,2,\dots$ Consequently, we have
\begin{align}\label{Case1vnh}
v^{(h)}(\cdot,t_n)-v^{(h)}_n =0 .
\end{align}

In Case 2', by using the Laplace transform, we can derive the following error representation (see Appendix \ref{AppendixB}):
\begin{align}\label{vhtn-vhn-}
v^{(h)}(\cdot,t_n)-v_n^{(h)} 
&= \frac{1}{2\pi\mathrm{i}}\int_{\Gamma_{\theta,\kappa}}
        e^{zt_n}z^{-1} (z^\alpha - \Delta_h)^{-1} \Delta_hP_h\psi_0 \d z  \nonumber  \\
&\quad 
  - \frac{1}{2\pi\mathrm{i} }\int_{\Gamma_{\theta,\kappa}^{(\tau)} }e^{zt_{n}} e^{-z\tau}\delta_\tau(e^{-z\tau})^{-1} (\delta_\tau(e^{-z\tau})^\alpha - \Delta_h)^{-1} \Delta_h P_h\psi_0   \,\d z \nonumber\\
&=
\frac{1}{2\pi\mathrm{i}} \int_{\Gamma_{\theta,\kappa}^{(\tau)} } e^{zt_n}
D_h^{(1)}(z)\Delta_h^\gamma P_h\psi_0 \,\d z \nonumber\\
&\quad
+\frac{1}{2\pi\mathrm{i}} \int_{\Gamma_{\theta,\kappa}\backslash\Gamma_{\theta,\kappa}^{(\tau)} } e^{zt_n} z^{ -1}\Delta_h^{1-\gamma}(z^\alpha - \Delta_h)^{-1}  \Delta_h^\gamma P_h\psi_0  \d z \nonumber\\
&=:
\mathcal{I}_n^{(h)} + \mathcal{J}_n^{(h)},
\end{align}
where $\gamma\in[0,1]$ and
\begin{align*}
&D_h^{(1)}(z)=z^{-1}\Delta_h^{1-\gamma}(z^\alpha - \Delta_h)^{-1} -e^{-z\tau}\delta_\tau(e^{-z\tau})^{-1} \Delta_h^{1-\gamma}(\delta_\tau(e^{-z\tau})^\alpha - \Delta_h)^{-1} .
\end{align*}
Lemma \ref{ineq-1} and Lemma \ref{Lemma:resolvent} imply
\begin{align*}
\| D_h^{(1)}(z)\|
&\le
C|z|^{-1}\big(\|\Delta_h^{1-\gamma}(z^\alpha-\Delta_h)^{-1}\|
+\|\Delta_h^{1-\gamma}(\delta_\tau(e^{-z\tau})^\alpha-\Delta_h)^{-1}\|\big) \nonumber\\\
&\le
C|z|^{-\gamma\alpha-1} \quad(\textrm{here we use }\eqref{resolv-frac-D})
\end{align*}and
\begin{align*}
\| D_h^{(1)}(z)\|
&\le
C|z|^{-1}\|\Delta_h^{1-\gamma} [(z^\alpha-\Delta_h)^{-1}-(\delta_\tau(e^{-z\tau})^\alpha-\Delta_h)^{-1}]\| 
\nonumber\\
&\quad
+C|z^{-1}-e^{-z\tau}\delta_\tau(e^{-z\tau})^{-1}| \|\Delta_h^{1-\gamma}(\delta_\tau(e^{-z\tau})^\alpha-\Delta_h)^{-1}\| 
\nonumber\\
&\le
C|z|^{-1} |z^\alpha-\delta_\tau(e^{-z\tau})^\alpha| \|\Delta_h^{1-\gamma}(\delta_\tau(e^{-z\tau})^\alpha-\Delta_h)^{-1}\| \|(z^\alpha-\Delta_h)^{-1}\| 
\nonumber\\
&\quad
+C|
(\delta_\tau(e^{-z\tau})-z)+(1-e^{-z\tau})z 
| |z|^{-1} |\delta_\tau(e^{-z\tau})|^{-1} |z|^{-\gamma\alpha} \quad (\textrm{use }\eqref{resolv-frac-D})
\nonumber\\
&\le
C|z|^{-\gamma\alpha}\tau.
\end{align*}
The last two inequalities further imply 
\begin{align*}
&\| D_h^{(1)}(z)\| \le C|z|^{-\gamma\alpha-1+\eta}  \tau^{\eta} \quad\forall\, \eta\in[0,1]. 
\end{align*}
Then, we have 
\begin{align}\label{Inhh}
 \| \mathcal{I}_n^{(h)}  \|
&\le C \int_{\Gamma_{\theta,\kappa}^{(\tau)} } |e^{zt_n}| 
\|D_h^{(1)}(z)\| \|\Delta_h^\gamma  P_h\psi_0 \|  \,|\d z| \nonumber \\
&\le C  \int_{\Gamma_{\theta,\kappa}^{(\tau)} } |e^{zt_n}| |z|^{-\gamma\alpha-1+\eta}  \tau^{\eta}
\| \Delta_h^\gamma P_h\psi_0\|  \,|\d z| \nonumber \\
&\le C \tau^{\eta} \| \Delta_h^\gamma  P_h\psi_0 \| 
\bigg(\int_{\kappa}^{\frac{\pi}{\tau\sin(\theta)}}  e^{rt_n\cos(\theta)}  
r^{-\gamma\alpha-1+\eta}  \,\d r   
+ \int_{-\theta}^{\theta} e^{\kappa t_n\cos(\varphi)}  
 \, \kappa^{-\gamma\alpha+\eta} \d \varphi \bigg) \nonumber \\
 &\le C \tau^{\eta}\| \Delta_h^\gamma  P_h\psi_0 \| 
\bigg( t_n^{\gamma\alpha-\eta} \int_{\kappa t_n}^{\frac{t_n \pi}{\tau\sin(\theta)}} e^{-s|\cos(\theta)|} 
s^{-\gamma\alpha-1+\eta}  \,\d s
 + \kappa^{-\gamma\alpha+\eta}  \bigg) \nonumber \\
 &\le C  \tau^{\gamma\alpha}(1+\kappa^{-\gamma\alpha+\eta}) \| \Delta_h^\gamma  P_h\psi_0\|    \,  \qquad(\mbox{this requires}\,\,\, \eta>\gamma\alpha) \nonumber \\ 
&\le C\tau^{\gamma\alpha} \|\psi_0\|_{\dot H^{2\gamma}(\mathcal{O})} ,  
\end{align}
where we have used \eqref{StabPh5} in the last inequality. 
Similarly, by using Lemma \ref{Lemma:resolvent} we have 
\begin{align}\label{Jnhh}
\| \mathcal{J}_n^{(h)}  \|
&\le C \int_{\Gamma_{\theta,\kappa}\backslash\Gamma_{\theta,\kappa}^{(\tau)} } 
|e^{zt_n} | |z|^{-1} \|\Delta_h^{1-\gamma}(z^\alpha - \Delta_h)^{-1}\|
 \|\Delta_h^\gamma  P_h\psi_0\| |\d z| \nonumber \\
&\le C \|\Delta_h^\gamma P_h\psi_0\|\int_{\frac{\pi}{\tau\sin(\theta)}}^\infty e^{-rt_n|\cos(\theta)|} 
r^{-1-\gamma\alpha}  \d r \quad (\textrm{here we use }\eqref{resolv-frac-D})
\nonumber  \\
&\le 
C\tau^{1+\gamma\alpha}\|\Delta_h^\gamma P_h\psi_0\| \int_{\frac{\pi}{\tau\sin(\theta)}}^\infty e^{-rt_n|\cos(\theta)|} \d r \nonumber \\
&\le 
C\tau^{\gamma\alpha} \|\Delta_h^\gamma P_h\psi_0\|  \nonumber \\
&\le C\tau^{\gamma\alpha} \|\psi_0\|_{\dot H^{2\gamma}(\mathcal{O})},
\end{align}
where we have used \eqref{StabPh5} again in the last inequality. 
By choosing 
$\gamma=\frac{1}{2\alpha}-\frac{d}{4}\in(0,1)$ and 
$\eta>\frac{1}{2}-\frac{\alpha d}{4}\in(0,1)$ 
(recall that $1<\alpha<\frac{2}{d}$), 
substituting \eqref{Inhh} and \eqref{Jnhh} into \eqref{vhtn-vhn-} yields 
\begin{align}\label{Case2vnh}
\|v^{(h)}(\cdot,t_n)-v_n^{(h)} \|  
&\le 
C\tau^{\frac{1}{2}-\frac{\alpha d}{4}} \|\psi_0\|_{\dot H^{\frac{1}{\alpha}-\frac{d}{2}}(\mathcal{O})}
\le C\tau^{\frac{1}{2}-\frac{\alpha d}{4}} \|\psi_0\|_{\dot H^{\chi}(\mathcal{O})}  ,
\end{align}
where we have noted that $\frac{1}{\alpha}-\frac{d}{2}\le \chi$ for $\alpha>1$.

In Case 3', with $\alpha\in(0,2)$, we have 
\begin{align} 
v^{(h)} (\cdot,t_n)
&=
 \int_{0}^{t_n} E^{(h)}(t_n-s) P_h  f(\cdot,s) \d s  ,  \label{vhsemi-expr1}\\
v_n^{(h)}  
&=
 \int_{0}^{t_n} E_\tau^{(h)}(t_n-s)P_h  f_\tau (\cdot,s) \d s  , \label{vhsemi-expr2}
\end{align} 
where $E_\tau^{(h)}(t) $ is given by \eqref{Etauht-def} and $f_\tau$ is the piecewise constant (in time) function defined in \eqref{f-tau}. 
The difference of the two expressions \eqref{vhsemi-expr1} and \eqref{vhsemi-expr2} yields  
\begin{align}\label{KnLn-}
v^{(h)} (\cdot,t_n)-v_n^{(h)}  
&=
 \int_{0}^{t_n} \big[E^{(h)}(t_n-s)-E_\tau^{(h)}(t_n-s)\big] P_hf_\tau (\cdot,s) \d s \nonumber \\
&\quad  +\int_{0}^{t_n} E^{(h)}(t_n-s) \big[P_h f(\cdot,s)-P_h f_\tau (\cdot,s)\big] \d s  \nonumber \\
&=: {\cal K}_n^{(h)} +{\cal L}_n^{(h)} .
\end{align} 
In the following, we estimate ${\cal K}_n^{(h)}$ and ${\cal L}_n^{(h)}$ separately. 

By the inverse Laplace transform rule $\mathcal{L}^{-1}(\widehat f \widehat g)(t)=\int_0^t \mathcal{L}^{-1}(\widehat f)(t-s)\mathcal{L}^{-1}(\widehat g)(s)\d s$, we have
\begin{align}\label{knhhh}
{\cal K}_n^{(h)}
&= \frac{1}{2\pi\mathrm{i}}  \int_{\Gamma_{\theta,\kappa}}
        e^{zt_n}z^{\alpha-1}(z^\alpha - \Delta_h)^{-1} P_h\widehat f_\tau(\cdot,z) \d z 
\nonumber \\
  &\quad 
  - \frac{1}{2\pi\mathrm{i} }\int_{\Gamma_{\theta,\kappa}^{(\tau)} }e^{zt_{n}} \frac{z\tau}{e^{z\tau}-1}\delta_\tau(e^{-z\tau})^{\alpha -1} \delta_\tau(e^{-z\tau})^\alpha - \Delta_h)^{-1} 
  P_h\widehat f_\tau (\cdot,z)  \,\d z 
\nonumber \\
&= \frac{1}{2\pi\mathrm{i}}  \int_{\Gamma_{\theta,\kappa}\backslash \Gamma_{\theta,\kappa}^{(\tau)}}
        e^{zt_n}  \widehat D_h^{(2)}(z) P_h\widehat f_\tau(\cdot,z) \d z 
  - \frac{1}{2\pi\mathrm{i} }\int_{\Gamma_{\theta,\kappa}^{(\tau)} }e^{zt_{n}}  \widehat D_h^{(3)}(z) 
  P_h\widehat f_\tau (\cdot,z)  \,\d z 
\nonumber\\
&=\int_0^{t_n} D_h^{(2)}(t_n-s) P_h f_\tau(\cdot,s) \d s 
+\int_0^{t_n} D_h^{(3)}(t_n-s) P_h f_\tau(\cdot,s) \d s ,
\end{align}
where 
\begin{align*}
& \widehat D_h^{(2)}(z) 
= z^{\alpha-1}(z^\alpha - \Delta_h)^{-1}  ,\\
& \widehat D_h^{(3)}(z) 
=z^{\alpha-1}(z^\alpha - \Delta_h)^{-1} 
- \frac{z\tau}{e^{z\tau}-1}\delta_\tau(e^{-z\tau})^{\alpha -1} \delta_\tau(e^{-z\tau})^\alpha - \Delta_h)^{-1} ,
\end{align*}and
\begin{align*}
D_h^{(2)}(t)\phi
=
\frac{1}{2\pi\ii}\int_{\Gamma_{\theta,\kappa}\backslash \Gamma_{\theta,\kappa}^{(\tau)}} e^{zt} \widehat D_h^{(2)}(z)\phi\d z ,\quad
D_h^{(3)}(t)\phi
=
\frac{1}{2\pi\ii}\int_{\Gamma_{\theta,\kappa}^{(\tau)}} e^{zt} \widehat D_h^{(3)}(z)\phi\d z .
\end{align*}
Using the similar method as introduced in \cite[Lemma 3.4]{GunzburgerLiWang2017} (as well as the inequalities \eqref{resolv-est-discr} in Lemma \ref{Lemma:resolvent} and \eqref{1-expztau-1} in Appendix \ref{AppendC}), it is easy to see that
\begin{align}\label{Dh2Dh3}
&\| \widehat D_h^{(2)}(z) \| +\| \widehat D_h^{(3)}(z) \|  \le C |z|^{-1} 
\qquad\mbox{and}\qquad 
\| \widehat D_h^{(3)}(z) \| \le C  \tau .
\end{align}
The last two inequalities further imply 
\begin{align}\label{est-Dh3}
&\|\widehat D_h^{(3)}(z) \| \le C |z|^{-\theta} \tau^{1-\theta} 
\quad\forall\,\theta\in[0,1]. 
\end{align}
Consequently, 
\begin{align}
\begin{aligned}
\|D_h^{(2)}(t) \| 
&\le C\int_{\Gamma_{\theta,\kappa}\backslash \Gamma_{\theta,\kappa}^{(\tau)}}
        |e^{zt}| \|  \widehat D_h^{(2)}(z)  \| | \d z |  \\
&\le C\int_{\frac{\pi}{\tau\sin(\theta)}}^\infty 
       e^{rt\cos(\theta)}  r^{-1} \d r  \\
&\le C\tau \int_{\frac{\pi}{\tau\sin(\theta)}}^\infty 
       e^{rt\cos(\theta)}  \d r  \\
&\le Ct^{-1}\tau\int_{\frac{t\pi}{\tau\sin(\theta)}}^\infty 
       e^{s\cos(\theta)}  \d s  \\
&\le Ct^{-1} \tau  \\
&\le Ct^{-1/q} \tau^{1/q} \qquad\forall\, q>1  , \,\,\, t\ge \tau ,
\end{aligned}
\end{align}
and 
\begin{align}
\begin{aligned}
\|D_h^{(2)}(t) \| 
&\le C\int_{\Gamma_{\theta,\kappa}\backslash \Gamma_{\theta,\kappa}^{(\tau)}}
        |e^{zt}| \|  \widehat D_h^{(2)}(z)  \| | \d z |  \\
&\le C\int_{\frac{\pi}{\tau\sin(\theta)}}^\infty 
       e^{rt\cos(\theta)}  r^{-1} \d r  \\
&\le C\int_{\frac{t\pi}{\tau\sin(\theta)}}^\infty 
       e^{s\cos(\theta)} s^{-1} \d s  \\
&\le C\int_1^\infty
       e^{s\cos(\theta)} s^{-1} \d s +C\int_{\frac{t\pi}{\tau\sin(\theta)}}^1
       e^{s\cos(\theta)} s^{-1} \d s  \\
&\le C+C\ln\big(t^{-1}\tau\big) \\
&\le Ct^{-1/q}\tau^{1/q} \qquad\forall\, q>1  , \,\,\, t\in(0, \tau) .
\end{aligned}
\end{align}
A combination of the last two estimates gives the following estimate of $\|D_h^{(2)}(t) \| $: 
\begin{align}\label{Est-Dh2}
\begin{aligned}
\|D_h^{(2)}(t) \| 
&\le Ct^{-1/q} \tau^{1/q} \qquad\forall\, q>1  , \,\,\, t\ge 0 .
\end{aligned}
\end{align}
Similarly, we have 
\begin{align}\label{Est-Dh3}
\begin{aligned}
\|D_h^{(3)}(t) \| 
&\le C\int_{\Gamma_{\theta,\kappa}^{(\tau)}}
        |e^{zt}| \| \widehat D_h^{(3)}(z)  \| | \d z |  \\
&\le C\tau^{1/q} \int_{\Gamma_{\theta,\kappa}^{(\tau)}}
        |e^{zt}| |z|^{-(1-1/q)} | \d z | \qquad\mbox{(set $\theta=1-1/q$ in \eqref{est-Dh3})} \\
&\le C\tau^{1/q} \int_\kappa^{\frac{\pi}{\tau\sin(\theta)}}  
       e^{rt\cos(\theta)} r^{-(1-1/q)} \d r 
       +C\tau^{1/q}  \int_{-\theta}^{\theta} 
       e^{\kappa t\cos(\varphi)}  \kappa^{1/q}\d \varphi  \\
&\le Ct^{-1/q}\tau^{1/q} \int_{\kappa t}^{\frac{t\pi}{\tau\sin(\theta)}}
       e^{s\cos(\theta)}   \d s +C \kappa^{1/q}\tau^{1/q} e^{\kappa T}  \\
&\le Ct^{-1/q}\tau^{1/q}+CT^{-1/q} \tau^{1/q}  \\
&\le Ct^{-1/q} \tau^{1/q} \qquad\forall\, q>1  , \,\,\, t\in(0,T] .
\end{aligned}
\end{align}
Substituting \eqref{Est-Dh2}-\eqref{Est-Dh3} into \eqref{knhhh} yields 
\begin{align} \label{Est-Knh-f}
\|{\cal K}_n^{(h)}\|
&=\bigg\|\int_0^{t_n} D^{(2)}(t_n-s) P_h f_\tau(\cdot,s) \d s 
+\int_0^{t_n} D^{(3)}(t_n-s) P_h f_\tau(\cdot,s) \d s \bigg\|  \nonumber \\
&\le C\tau^{1/q} 
\int_0^{t_n}(t_n-s)^{-1/q}\|P_h f_\tau(\cdot,s)\| \d s \nonumber \\
&\le C\tau^{1/q}\|f_\tau\|_{L^{q',1}(0,t_n;L^2(\mathcal{O}))}   \nonumber \\
&\le C\tau^{1/q}\|f\|_{L^{q',1}(0,t_n;L^2(\mathcal{O}))}, 
\end{align}
where the last inequality follows from \eqref{Lp1-stable}. This completes the estimate of the first term in \eqref{KnLn-}.

The second term in \eqref{KnLn-} can be estimated as follows: 
\begin{align} \label{Est-Lnh-f}
\|{\cal L}_n^{(h)}\|
&=\bigg\| \sum_{j=1}^n\int_{t_{j-1}}^{t_j} E^{(h)}(t_n-s) \big[P_h f(\cdot,s)-P_h f_\tau (\cdot,s)\big] \d s \bigg\|  
 \nonumber \\
&=\bigg\| \sum_{j=1}^n\int_{t_{j-1}}^{t_j} 
\big(E^{(h)}(t_n-s)-E^{(h)}(t_n-t_{j-1})\big) \big[P_h f(\cdot,s)-P_h f_\tau (\cdot,s)\big] \d s \bigg\|  \nonumber  \\
&\le \sum_{j=1}^n\int_{t_{j-1}}^{t_j} 
\|E^{(h)}(t_n-s)-E^{(h)}(t_n-t_{j-1})\| 
\|P_h f(\cdot,s)-P_h f_\tau (\cdot,s)\| \d s   \nonumber \\
&\le \sum_{j=1}^n\int_{t_{j-1}}^{t_j} 
C\tau^{1/q}(t_n-s)^{-1/q}  
\|P_h f(\cdot,s)-P_h f_\tau (\cdot,s)\| \d s  \nonumber  \\
&\le C\tau^{1/q}\|f-f_\tau\|_{L^{q',1}(0,t_n;L^2(\mathcal{O}))}   \nonumber \\
&\le C\tau^{1/q}\|f\|_{L^{q',1}(0,t_n;L^2(\mathcal{O}))} ,
\end{align}
where we have used, for $q>1$, 
\begin{align}\label{Ehtns-j}
&
\|E^{(h)}(t_n-s)-E^{(h)}(t_n-t_{j-1})\| \nonumber \\
&\le
C \int_{\Gamma_{\theta,\kappa}} |e^{z(t_n-s)}-e^{z(t_n-t_{j-1})}| |z|^{\alpha-1}\|(z^\alpha-\Delta_h)^{-1}\| |\d z|  \nonumber \\
&\le
C \int_{\Gamma_{\theta,\kappa}} |e^{z(t_n-s)}| |1-e^{z(s-t_{j-1})}| |z|^{-1} |\d z|  
\qquad(\textrm{here we use }\eqref{resolv-est-discr})
\nonumber \\
&\le
C \int_{\Gamma_{\theta,\kappa}} |e^{z(t_n-s)}| \tau^{1/q}|z|^{-(1-1/q)} |\d z| 
\qquad\mbox{(here we use \eqref{1-expztau-2})} \nonumber  \\
&\le
C\tau^{1/q}(t_n-s)^{-1/q}.
\end{align}
By choosing $q>1$ to satisfy $\frac{1}{q}=\frac{1}{2}-\frac{\alpha d}{4}$, we have $q'=\frac{4}{2+\alpha d}$. 
Then substituting \eqref{Est-Knh-f} and \eqref{Est-Lnh-f} into \eqref{KnLn-} yields 
\begin{align}\label{Case3vnh}
\|v_n^{(h)}-v^{(h)} (\cdot,t_n)  \|
&\le C\tau^{\frac{1}{2}-\frac{\alpha d}{4}}\|f\|_{L^{\frac{4}{2+\alpha d},1}(0,t_n;L^2(\mathcal{O}))}   ,
\end{align}
completing the proof of Case 3. 

By combining \eqref{Case1vnh}, \eqref{Case2vnh} and \eqref{Case3vnh} (the results of Cases 1, 2 and 3), we obtain 
\begin{align}\label{vnvvvv}
\|v_n^{(h)}-v^{(h)} (\cdot,t_n)\| 
&\le 
C\big(\|\psi_0\|_{\dot H^{\chi}(\mathcal{O})} +\|f\|_{L^{\frac{4}{2+\alpha d},1}(0,t_n;L^2(\mathcal{O}))}\big) \tau^{\frac{1}{2}-\frac{\alpha d}{4}}  .
\end{align}
The estimates \eqref{vvvvvv} and \eqref{vnvvvv} imply \eqref{error-v-vnh}, completing the proof of Theorem \ref{MainTHM}.\qed

\section{Numerical examples}
In this section, we present numerical examples to illustrate the theoretical analyses. 

We consider the one-dimensional stochastic partial integro-differential equation \eqref {Frac-SPDE} for $0\le x\le 1,\ 0<t\le 1$, with homogeneous Dirichlet boundary condition and initial condition $\psi_0(x)=x(1-x)$. Here, we let $\sigma=1$ in \eqref{Frac-SPDE} and 
\begin{align*}
f(x,t)
=
\left\{\begin{aligned}
&1
&&\mbox{for }0\le x\le\frac{1}{2},\\
&-1
&&\mbox{for }\frac{1}{2}<x\le 1.
\end{aligned}\right.  
\end{align*}
The problem \eqref{Frac-SPDE} is discretized by using the scheme \eqref{fully-discrete}-\eqref{fully-discrete-2}. 

To investigate the convergence rate in space, we first solve the problem \eqref{Frac-SPDE} by taking the mesh size $h_k=1/M_k=2^{-k}$, $k=2,3,4,5$, and using a sufficiently small time step $\tau=2^{-14}$ so that the temporal discretization error is relatively negligible. Then, the error 
\begin{align}
E(h_k)=\frac{1}{I}\sum_{i=1}^I \|\psi_N^{(h_k)}(\omega_i)-\psi_N^{(h_{k-1})}(\omega_i)\|
\end{align}is computed for $k=3,4,5$, by using $I=10000$ independent realizations for each spatial mesh size. 
By Theorem \ref{MainTHM}, the error $E(h_k)$ is expected to have the convergence rate $O(h^{\frac{1}{\alpha}-\frac{1}{2}})$ for $\alpha\in[\frac{1}{2},2)$, and $O(h^{\frac{3}{2}})$ for $\alpha\in(0,\frac{1}{2})$ in one-dimensional spatial domain. 
The numerical results are presented in Table \ref{space-error} and consistent with the theoretical analyses that the spatial order of convergence increases as $\alpha$ decreases and stays at the order $\frac{3}{2}$ when $\alpha$ further decreases. 

Secondly, we solve the problem \eqref{Frac-SPDE} by using the time step $\tau_k=2^{-k}$, $k=6,7,8,9$. In order to focus on the temporal discretization error, a sufficiently small spatial mesh size $h=1/M=2^{-10}$ is used such that the spatial discretization error can be relatively negligible.  Similarly, we consider $I=10000$ independent realizations for each time step and compute the error $E(\tau_k)$ by  
\begin{align}
E(\tau_k)=\frac{1}{I}\sum_{i=1}^I \|\psi_{N,\tau_k}^{(h)}(\omega_i)-\psi_{N,\tau_{k-1}}^{(h)}(\omega_i)\|
\end{align}for $k=7,8,9$, 
which is expected to have the convergence rate $O(\tau^{\frac{1}{2}-\frac{\alpha}{4}})$ by Theorem \ref{MainTHM}. Clearly, the results in Table \ref{time-error} illustrate the sharp convergence rate.


\begin{table}[!ht]
\begin{center}
\caption{$E(h_k)$ and convergence rates in space }
\label{space-error}
\begin{tabular}{c|ccc|c}
\hline
\hline
\ {$\alpha\backslash h_k$}&$2^{-3}$&$2^{-4}$&$2^{-5}$&order\  \\
\hline
\ $\alpha=0.25$ &1.1669e-02 &3.9124e-03 &1.3519e-03 &1.555 (1.500)\ \\
\ $\alpha=0.75$ &2.4353e-02 &1.2987e-02 &6.6322e-03  &0.938 (0.833)\ \\
\ $\alpha=1.25$ &8.3694e-02 &6.7186e-02 &5.4196e-02 &0.314 (0.300)\ \\
\hline
\hline
\end{tabular} 
\end{center}
\end{table}

\begin{table}[!ht]
\begin{center}
\caption{$E(\tau_k)$ and convergence rates in time }
\label{time-error}
\begin{tabular}{c|ccc|c}
\hline
\hline
\ {$\alpha\backslash \tau_k$}&$2^{-7}$&$2^{-8}$&$2^{-9}$&order\  \\
\hline
\ $\alpha=0.25$ &2.2103e-03 &1.7275e-03 &1.3454e-03 &0.359 (0.4375)\ \\
\ $\alpha=0.75$ &1.5613e-02 &1.2621e-02 &1.0177e-02  &0.309 (0.3125)\ \\
\ $\alpha=1.25$ &5.0056e-02 &4.4012e-02 &3.8869e-02 &0.183 (0.1875)\ \\
\hline
\hline
\end{tabular} 
\end{center}
\end{table}



\appendix

\section{Mild solution of (\ref{Frac-SPDE})}\label{AppendixA}

In the case $\alpha\in(1,2)$, the boundary condition $\partial_t^{1-\alpha}\psi=0$ is equivalent to $\psi=0$ on $\partial\Omega$ (this can be checked by taking Laplace transform in time). 
Similarly, in the case $\alpha\in(0,1]$, the boundary condition $\partial_t^{1-\alpha}\psi=0$ is equivalent to $\psi-\psi_0=0$ on $\partial\Omega\times[0,\infty)$, where $\psi_0=\psi(\cdot,0)$ is the initial value in \eqref{Frac-SPDE}. 

In the case $\sigma=0$, the solution of the corresponding deterministic problem of \eqref{Frac-SPDE} can be expressed by (via Laplace transform, cf. \cite[(3.11) and line 4 of page 12]{LubichSloanThomee:1996} in the case $\alpha\in(1,2)$) 
\begin{align}\label{Deter-sol-repr}
\psi(\cdot,t)
&= 
\left\{\begin{aligned}
&\psi_0 + \int_0^t  E(t-s) f(\cdot,s)\d s            &&\mbox{if}\,\,\,\alpha\in(0,1] , \\
&E(t) \psi_0 + \int_0^t  E(t-s) f(\cdot,s)\d s     &&\mbox{if}\,\,\,\alpha\in(1,2) ,  
\end{aligned}
\right.
\end{align}
where the operator $E(t):L^2(\mathcal{O})\rightarrow L^2(\mathcal{O})$ is given by 
\begin{equation}\label{eqn:EF}
E(t) \phi:=\frac{1}{2\pi {\rm i}}\int_{\Gamma_{\theta,\kappa}}e^{zt} z^{\alpha-1} (z^\alpha-\Delta )^{-1}\phi\, \d z  \quad
\forall\, \phi\in L^2(\mathcal{O}) 
\end{equation}
with integration over a contour $\Gamma_{\theta,\kappa}$ on the complex plane. 

Correspondingly, the mild solution of the stochastic problem \eqref{Frac-SPDE} is defined as 
(cf. \cite[Proposition 2.7]{KP} and \cite{MijenaNane2015})
\begin{align}\label{Mild-sol-}
\psi(\cdot,t)
&= 
\left\{\begin{aligned}
&\psi_0 + \int_0^t  E(t-s) f(\cdot,s)\d s   +\sigma\int_0^t  E(t-s)\d W(\cdot,s)         &&\mbox{if}\,\,\,\alpha\in(0,1] , \\
&E(t) \psi_0 + \int_0^t  E(t-s) f(\cdot,s)\d s  +\sigma\int_0^t  E(t-s)\d W(\cdot,s)   &&\mbox{if}\,\,\,\alpha\in(1,2) .  
\end{aligned}
\right.
\end{align}
For any given initial data $\psi_0\in L^2(\mathcal{O})$ and source $f\in L^1(0,T;L^2(\mathcal{O}))$,  
the expression \eqref{Mild-sol-} defines a mild solution $\psi\in C([0,T];L^2(\Omega; L^2(\mathcal{O})))$.  
In the case $\psi_0=f=0$ and $\sigma\neq 0$, a simple proof of this result can be found in \cite[Appendix]{GunzburgerLiWang2017}; in the case $\sigma=0$ ($\psi_0$ and $f$ may not be zero), the result is a consequence of the boundedness of the operator $E(t):L^2(\mathcal{O})\rightarrow L^2(\mathcal{O})$, i.e., 
\begin{align}\label{BD-Et}
\|E(t)v\|
&\le
C\int_{\Gamma_{\theta,\kappa}} |e^{zt}| |z|^{\alpha-1} \|(z^\alpha-\Delta)^{-1}v\| |\d z| \nonumber \\
&\le
C\|v\|\int_\kappa^{\infty} e^{-rt|\cos(\theta)|} r^{-1} \d r
+C\|v\|\int_{-\theta}^\theta e^{\kappa t\cos(\varphi)} \d\varphi \nonumber \\
&\le C\|v\| \qquad \forall\, v\in L^2(\mathcal{O}).
\end{align}
Similarly, the discrete operator $E^{(h)}(t):X_h\rightarrow X_h$ defined by 
\begin{equation}\label{eqn:Eht}
E^{(h)}(t) \phi:=\frac{1}{2\pi {\rm i}}\int_{\Gamma_{\theta,\kappa}}e^{zt} z^{\alpha-1} (z^\alpha-\Delta_h )^{-1}\phi\, \d z   \quad
\forall\, \phi\in X_h ,
\end{equation}
is also bounded on the finite element subspace $X_h$, i.e.,  
\begin{align}\label{BD-Eht}
\|E^{(h)}(t)v\|\le C\|v\| \quad\forall\, v\in X_h ,
\end{align}
where the constant $C$ is independent of the mesh size $h$.

\section{Representation of the discrete solutions}\label{AppendixB}

For $f=0$ we prove the following representation of the solutions of \eqref{semi-FEM} and \eqref{Deter-discreteh-2} : 
\begin{align}
&v^{(h)}(\cdot,t_n)
= P_h\psi_0 + \frac{1}{2\pi\mathrm{i}}\int_{\Gamma_{\theta,\kappa}}
        e^{zt_n}z^{-1} (z^\alpha - \Delta_h)^{-1} \Delta_hP_h\psi_0 \d z ,
        \label{Repr-vh} \\
&v_n^{(h)} 
=P_h\psi_0 + \frac{1}{2\pi i}\int_{\Gamma_{\theta,\kappa}^{(\tau)}} e^{t_nz}e^{-z\tau} \delta(e^{-z\tau})^{-1} \big(\delta(e^{-z\tau})^{\alpha}-\Delta_h\big)^{-1} \Delta_h P_h\psi_0   \,  \d  z  ,
\label{Tyounrepr2}
\end{align}
which are used in \eqref{vhtn-vhn-} in estimating the error of temporal discretization. 

In fact, \eqref{Repr-vh} is a consequence of \eqref{Mild-sol-}: replacing $E(t)$ by $E^{(h)}(t)$ and substituting $\phi=P_h\psi_0$ yield 
\begin{align*}
\begin{aligned}
v^{(h)}(\cdot,t_n)
&=\frac{1}{2\pi {\rm i}}\int_{\Gamma_{\theta,\kappa}}e^{zt} z^{\alpha-1} (z^\alpha-\Delta_h )^{-1}P_h\psi_0\, \d z \\
&=\frac{1}{2\pi {\rm i}}\int_{\Gamma_{\theta,\kappa}}e^{zt} z^{-1}(z^\alpha-\Delta_h+\Delta_h) (z^\alpha-\Delta_h )^{-1}P_h\psi_0\, \d z \\
&=\frac{1}{2\pi {\rm i}}\int_{\Gamma_{\theta,\kappa}}e^{zt} z^{-1} P_h\psi_0\, \d z  + \frac{1}{2\pi {\rm i}}\int_{\Gamma_{\theta,\kappa}}e^{zt} z^{-1}(z^\alpha-\Delta_h )^{-1}\Delta_hP_h\psi_0\, \d z \\
&=P_h\psi_0 + \frac{1}{2\pi {\rm i}}\int_{\Gamma_{\theta,\kappa}}e^{zt} z^{-1}(z^\alpha-\Delta_h )^{-1}\Delta_hP_h\psi_0\, \d z ,
\end{aligned} 
\end{align*}
where we have used the identity $\frac{1}{2\pi {\rm i}}\int_{\Gamma_{\theta,\kappa}}e^{zt} z^{-1}  \d z =1$ (i.e., the inverse Laplace transform of $z^{-1}$ is $1$). 

It remains to prove \eqref{Tyounrepr2}. To this end, we rewrite \eqref{Deter-discreteh-2} as 
\begin{equation}\label{Rewrite-vn_h0}
\bar\partial_\tau (v_n^{(h)}-P_h\psi_0)  -\Delta_h \bar\partial_\tau^{1-\alpha} (v_n^{(h)}-P_h\psi_0)
= \Delta_h  \bar\partial_\tau^{1-\alpha} (P_h\psi_0)_n ,
\end{equation} 
where $\bar\partial_\tau^{1-\alpha} (P_h\psi_0)_n:=\frac{1}{\tau^{1-\alpha}}\sum_{j=1}^n b_{n-j} P_h\psi_0$. 
Since we are only interested in the solutions $v_n^{(h)}$, $n=1,\dots,N$, we define 
$$
\widetilde v_n^{(h)}
=\left\{
\begin{array}{ll}
v_n^{(h)}   &\quad 1\le n\le N,\\[5pt]
P_h\psi_0 &\quad n\ge N+1 ,
\end{array}\right.
$$
which satisfies the equation 
\begin{equation}\label{Rewrite-vn_h}
\bar\partial_\tau (\widetilde v_n^{(h)}-P_h\psi_0)  -\Delta_h \bar\partial_\tau^{1-\alpha} (\widetilde v_n^{(h)}-P_h\psi_0)
= \Delta_h  \bar\partial_\tau^{1-\alpha} (P_h\psi_0)_n 
+ g_n , 
\end{equation}
with $g_n=0$ for $1\le n\le N$. The right-hand side of \eqref{Rewrite-vn_h} differs from \eqref{Rewrite-vn_h0} only for $n\ge N+1$, that 
\begin{align}
\|g_n\|
&\le
\|\Delta_h \bar\partial_\tau^{1-\alpha} (\widetilde v_n^{(h)}-P_h\psi_0)\|
+\|\Delta_h  \bar\partial_\tau^{1-\alpha} (P_h\psi_0)_n \| \nonumber\\
&\le
\frac{1}{\tau^{1-\alpha}}\sum_{j=1}^N |b_{n-j}|\|\widetilde v_j^{(h)}-P_h\psi_0\|
+\frac{1}{\tau^{1-\alpha}}\sum_{j=1}^n |b_{n-j}|\|P_h\psi_0\| \nonumber\\
&\le
C\tau^{\alpha-1}\Big(\sum_{j=1}^N |b_{n-j}|+\sum_{j=1}^n |b_{n-j}|\Big) \nonumber\\
&\le
C\tau^{\alpha-1}n^{\alpha-1}, 
\end{align}
as $n\rightarrow \infty$. Thus $\sum_{n=N+1}^\infty g_n\zeta^n$ is an analytic function of $\zeta$ for $|\zeta|<1$. 
 
By \eqref{generate-dv}, summing up \eqref{Rewrite-vn_h} times $\zeta^n$ for $n=1,2,\dots$, yields 
\begin{equation}
\bigg(\frac{1-\zeta}{\tau}-\bigg(\frac{1-\zeta}{\tau}\bigg)^{1-\alpha}\Delta_h \bigg)
\sum_{n=1}^\infty(\widetilde v_n^{(h)}-P_h\psi_0)\zeta^n= \Delta_h  \bigg(\frac{1-\zeta}{\tau}\bigg)^{1-\alpha} \frac{\zeta}{1-\zeta}  P_h\psi_0 
+ \sum_{n=N+1}^\infty g_n\zeta^n  ,
\end{equation}
which implies 
\begin{align}
\begin{aligned}
\sum_{n=1}^\infty(\widetilde v_n^{(h)}-P_h\psi_0)\zeta^n
=
& \bigg(\frac{1-\zeta}{\tau}\bigg)^{-1} \bigg(\bigg(\frac{1-\zeta}{\tau}\bigg)^{\alpha}-\Delta_h \bigg)
^{-1}\Delta_h P_h\psi_0  \frac{\zeta}{\tau}  \\
& + \bigg(\frac{1-\zeta}{\tau}\bigg)^{\alpha-1}\bigg(\bigg(\frac{1-\zeta}{\tau}\bigg)^{\alpha}-\Delta_h \bigg)^{-1}\sum_{n=N+1}^\infty g_n\zeta^n .
\end{aligned}
\end{align}
For $\kappa>0$ and $\varrho_\kappa=e^{-(\kappa+1)\tau} \in(0,1)$, the Cauchy integral formula implies that 
\begin{align}\label{Tyounrepr0}
\begin{split}
&\widetilde v_n^{(h)}-P_h\psi_0 \\
&= \frac{1}{2\pi i}\int_{|\zeta|=\varrho_\kappa} \zeta^{-n-1}\sum_{n=1}^\infty (v_n^{(h)}-P_h\psi_0)\zeta^n  \d \zeta \\
&= \frac{1}{2\pi i}\int_{|\zeta|=\varrho_\kappa} \zeta^{-n} \bigg(\frac{1-\zeta}{\tau}\bigg)^{-1} \bigg(\bigg(\frac{1-\zeta}{\tau}\bigg)^{\alpha}-\Delta_h\bigg)^{-1}\frac{1}{\tau}\Delta_h P_h\psi_0    \d \zeta \\
&\quad + \frac{1}{2\pi i}\int_{|\zeta|=\varrho_\kappa} \bigg(\frac{1-\zeta}{\tau}\bigg)^{\alpha-1}\bigg(\bigg(\frac{1-\zeta}{\tau}\bigg)^{\alpha}-\Delta_h \bigg)^{-1}\sum_{m=N+1}^\infty g_m\zeta^{m-n-1}     \d \zeta . 
\end{split}
\end{align}
For $1\le n\le N$ the function $ \big(\frac{1-\zeta}{\tau}\big)^{\alpha-1}\big(\big(\frac{1-\zeta}{\tau}\big)^{\alpha}-\Delta_h \big)^{-1}\sum_{m=N+1}^\infty g_m\zeta^{m-n-1} $ is analytic in $|\zeta|<1$. Consequently, Cauchy's integral theorem implies 
$$
 \frac{1}{2\pi i}\int_{|\zeta|=\varrho_\kappa} \bigg(\frac{1-\zeta}{\tau}\bigg)^{\alpha-1}\bigg(\bigg(\frac{1-\zeta}{\tau}\bigg)^{\alpha}-\Delta_h \bigg)^{-1}\sum_{m=N+1}^\infty g_m\zeta^{m-n-1}     \d \zeta =0 .
$$
Substituting this identity into \eqref{Tyounrepr0} yields, for $1\le n\le N$, 
\begin{align}\label{Tyounrepr}
\begin{split}
&\widetilde v_n^{(h)}-P_h\psi_0\\
&=v_n^{(h)}-P_h\psi_0 \\
&= \frac{1}{2\pi i}\int_{|\zeta|=\varrho_\kappa} \zeta^{-n} \bigg(\frac{1-\zeta}{\tau}\bigg)^{-1} \bigg(\bigg(\frac{1-\zeta}{\tau}\bigg)^{\alpha}-\Delta_h\bigg)^{-1}\frac{1}{\tau}\Delta_h P_h\psi_0    \d \zeta \\
&= \frac{1}{2\pi i}\int_{\Gamma^\tau}e^{t_nz}e^{-z\tau} \bigg(\frac{1-e^{-\tau z}}{\tau}\bigg)^{-1} \bigg(\bigg(\frac{1-e^{-\tau z}}{\tau}\bigg)^{\alpha}-\Delta_h\bigg)^{-1} \Delta_h P_h\psi_0   \,  \d  z \\
&= \frac{1}{2\pi i}\int_{\Gamma^\tau}e^{t_nz}e^{-z\tau} \delta(e^{-z\tau})^{-1} \big(\delta(e^{-z\tau})^{\alpha}-\Delta_h\big)^{-1} \Delta_h P_h\psi_0   \,  \d  z ,
\end{split}
\end{align}
where we have used the change of variable $\zeta =e^{-z\tau}$, which converts the path of integration to the contour
\begin{align}\label{cGammatau}
\Gamma^\tau 
&=\left\{ z=\kappa+1+\mathrm{i} y: \, y\in{\mathbb R}\,\,\,\mbox{and}\,\,\,|y|\le {\pi}/{\tau} \right\} .
\end{align}
The angle condition \eqref{angle-delta} and \cite[Theorem 3.7.11]{ABHN} imply that the integrand on the right-hand side of \eqref{Tyounrepr} is analytic in the region 
\begin{align}
&\Sigma_{\theta,\kappa}^\tau = \Big \{z\in{\mathbb C} :
|{\rm arg}(z)|\le \theta,\,\,\, |z|\ge \kappa,\,\,|{\rm Im}(z)|\le \frac{\pi}{\tau} ,\,\,\,
{\rm Re}(z)\le \kappa +1  \Big \} ,
\label{Sigma-theta-k}
\end{align}
enclosed by the four paths $\Gamma^\tau$, $\Gamma_{\theta,\kappa}^{(\tau)}$ and ${\mathbb R}\pm \mathrm{i}\pi/\tau$, where $\Gamma_{\theta,\kappa}^{(\tau)}
\, =\left\{z\in \Gamma_{\theta,\kappa} : |{\rm Im}(z)|\le \frac{\pi}{\tau}\right\}$. 
Then Cauchy's theorem allows us to deform the integration path from $\Gamma^\tau$ to $\Gamma_{\theta,\kappa}^{(\tau)}$ in the integral \eqref{Tyounrepr} (the integrals on ${\mathbb R}\pm \mathrm{i}\pi/\tau$ cancels each other). This yields the desired representation \eqref{Tyounrepr2}.

\section{Some inequalities}\label{AppendC}

In this appendix, we prove the following two inequalities:
\begin{align}
&C_0^{\#}|z|\tau\le |1-e^{z\tau}|\le C_1^{\#}|z|\tau , &&\forall\, z\in \Gamma_{\theta,\kappa}^{(\tau)}, \label{1-expztau-1}\\
&|1-e^{z\tau}|\le C|z|^{1/q}\tau^{1/q}, &&\forall\, z\in\Gamma_{\theta,\kappa} , \,\,\,1\le q\le\infty .\label{1-expztau-2}
\end{align}
which have been used in \eqref{E_tautphij}, \eqref{Dh2Dh3} and \eqref{Ehtns-j}. 

{\bf Proof of \eqref{1-expztau-1}.} 
Note that 
\begin{align}
\Gamma_{\theta,\kappa}^{(\tau)}
&=\left\{z\in \mathbb{C}: |z|=\kappa  ,\, |\arg z|\le \theta\right\} \cup
  \left\{z\in \mathbb{C}: z=\rho e^{\pm {\rm i}\theta}, \rho\ge \kappa , |{\rm Im}(z)|\le \frac{\pi}{\tau} \right\}  \nonumber \\
&=: \Gamma_{\theta,\kappa}^{(\tau),1}\cup \Gamma_{\theta,\kappa}^{(\tau),2} .
\end{align}
For $z\in \Gamma_{\theta,\kappa}^{(\tau)}$ we have $|z|\tau\le \pi/\sin(\theta)$. Since $|z|\tau$ is bounded, the following Taylor expansion holds:
\begin{align}\label{1-eztau22}
&1-e^{z\tau} =-z\tau+O(|z|^2\tau^2) , 
\end{align}
which imiplies
\begin{align}
|1-e^{z\tau}| \le C_1^{+}|z|\tau ,\quad\mbox{if}\,\,\,z\in \Gamma_{\theta,\kappa}^{(\tau)}.
\end{align}
This proves the right-half inequality of \eqref{1-expztau-1}. 

From \eqref{1-eztau22} we also see that there exists a small constant $\gamma$ such that 
\begin{align}\label{1-eztau-small}
C_0^{+}|z|\tau\le |1-e^{z\tau}|  ,\quad\mbox{if}\,\,\,z\in \Gamma_{\theta,\kappa}^{(\tau)},\,\,\, |z|\tau<\gamma ,
\end{align}
If $|z|\tau\ge \gamma$, then the following inequality holds for $\theta$ satisfying the condition of Lemma \ref{ineq-1}:
$$
\gamma\le |z|\tau \le \frac{\pi}{\sin(\theta)}\le \pi \sqrt{1+4/\pi^2} \le \frac{3}{2}\pi .
$$
Since the function $g(w):=|1-e^{w}|$ is not zero for $\gamma\le |w| \le \frac{3}{2}\pi $, the function $g(w)$ must have a positive minimum value $\varpi$ for $\gamma\le |w| \le \frac{3}{2}\pi $, $i.e.$, $g(w)\ge \varpi$. Consequently, we have 
\begin{align}\label{1-eztau-big}
\begin{aligned}
&\varpi \frac{\sin(\theta)}{\pi}|z|\tau  
\le \varpi \le |1-e^{z\tau}| ,
&&\mbox{if}\,\,\,z\in \Gamma_{\theta,\kappa}^{(\tau)} ,\,\,\, |z|\tau\ge\gamma.  
\end{aligned}
\end{align}
where we have used the inequality $\frac{\sin(\theta)}{\pi}|z|\tau  \le 1$ in the last inequality. 
Combining \eqref{1-eztau-small} and \eqref{1-eztau-big} yields \eqref{1-expztau-1}.\qed 

{\bf Proof of \eqref{1-expztau-2}.} 
If $z\in\Gamma_{\theta,\kappa}$ and $|z|\tau\le \pi/\sin(\theta)$, then $z\in\Gamma_{\theta,\kappa}^{(\tau)}$. In this case, \eqref{1-expztau-1} implies 
\begin{align}
&|1-e^{z\tau}|\le C|z|\tau, &&\forall\, z\in\Gamma_{\theta,\kappa} , \,\,\, |z|\tau\le \pi/\sin(\theta),\\
&|1-e^{z\tau}|\le C, &&\forall\, z\in\Gamma_{\theta,\kappa} , \,\,\, |z|\tau\le \pi/\sin(\theta) .
\end{align}
The combination of the two inequalities above yields 
\begin{align}\label{expztau-2-small}
&|1-e^{z\tau}|\le C|z|^{1/q}\tau^{1/q}, &&\forall\, z\in\Gamma_{\theta,\kappa} , \,\,\, |z|\tau\le \pi/\sin(\theta).
\end{align}
If $z\in\Gamma_{\theta,\kappa}$ and $|z|\tau\ge \pi/\sin(\theta)$, then 
$$|e^{z\tau}|=e^{-|z|\tau\cos(\theta)}\le e^{-\pi/\tan(\theta)} ,$$
which implies 
\begin{align}\label{expztau-2-big}
&|1-e^{z\tau}|\le 1+e^{-\pi /\tan(\theta)}
\le 2 \le 2 \bigg(\frac{\sin(\theta)}{\pi}\bigg)^{\frac1q} |z|^{1/q}\tau^{1/q} ,
&&\forall\, z\in\Gamma_{\theta,\kappa} , \,\,\, |z|\tau\ge \pi/\sin(\theta).
\end{align}
Combining \eqref{expztau-2-small} and \eqref{expztau-2-big} yields \eqref{1-expztau-2}.\qed
 
\vspace{0.2in}
{\bf Funding.} 
The research of M. Gunzburger and J. Wang was supported in part by the USA National Science Foundation grant DMS-1315259, by the USA Air Force Office of Scientific Research grant FA9550-15-1-0001, and by a USA Defense Advanced Projects Agency contract administered under the Oak Ridge National Laboratory subcontract 4000145366. The research of B. Li was partially supported by the Hong Kong RGC grant 15300817.

\bibliographystyle{abbrv}
\bibliography{frac_SPDE}

\end{document}